\documentclass[reqno]{amsart}
\usepackage[margin = 1in]{geometry}
\usepackage{amsmath, amssymb, amsthm, fancyhdr, verbatim, graphicx}
\usepackage{enumerate}
\usepackage[all]{xy}
\usepackage[dvipsnames]{xcolor}
\usepackage{mathrsfs}
\usepackage{tikz-cd}
\usepackage{framed}
\usepackage[hidelinks, pagebackref]{hyperref}
\hypersetup{colorlinks=true, linkcolor=blue, citecolor=teal} 
\usepackage[titletoc]{appendix}
\usepackage{bbm}
\usepackage{lipsum}
\usepackage{adjustbox}
\usepackage{commands}
\usepackage{envi}
\usepackage{stmaryrd}
\usepackage{leftindex}
\usepackage[english]{babel}
\usepackage{dynkin-diagrams}
\usepackage[mathscr]{eucal}

\numberwithin{equation}{section}

\makeatletter
\def\th@remark{%
  \thm@headfont{\bfseries}%
  \normalfont 
  \thm@preskip \thm@preskip 
  \thm@postskip\thm@preskip
}
\def\imod#1{\allowbreak\mkern5mu({\operator@font mod}\,\,#1)}
\makeatother

\numberwithin{equation}{section}

\title{Gross--Prasad conjecture in geometric Langlands}
\author{Kaiyi Chen, Zeyu Wang}

\address{Massachusetts Institute of Technology, Department of Mathematics, 77 Massachusetts Avenue, Cambridge, MA 02139, USA}
\email{kychen03@mit.edu}
\address{Massachusetts Institute of Technology, Department of Mathematics, 77 Massachusetts Avenue, Cambridge, MA 02139, USA}
\email{wangzeyu@mit.edu}

\begin{document}

\begin{abstract}
We prove part of the global Gross--Prasad conjecture for $G=\SO_n \times \SO_{n-1}$ and $H=\SO_{n-1}$ in the unramified de Rham geometric Langlands setting proposed by Ben-Zvi--Sakellaridis--Venkatesh. 
\end{abstract}

\maketitle
\tableofcontents

\section{Introduction}\label{sec:intro}
\subsection{Setup}
\subsubsection{Automorphic side}
Let $C$ be a connected smooth projective curve over $\CC$. Let $H=\SO_{n-1}$ and $H'=\SO_n$ with $n\geq 4$, and let $G=H'\times H$. Let $X=G/H$ be the corresponding spherical $G$-variety. We use $T\sub B=NT\sub G$ to denote the subgroups as usual. Let \[\Bun_G:=\Map(C,\mathbb{B}G)\] be the moduli stack of $G$-bundles on $C$, and let \[\Bun_G^X:=\Bun_H=\Map(C,[X/G]).\] There is a natural map \[\pi_X:\Bun_G^X\to \Bun_G.\] 
Following \cite[\S10.3]{BZSV}, one can define the (unnormalized) \emph{period sheaf}
\begin{equation}\label{eq:periodsheaf}
\cP_X := (\pi_{X,\dR})_!\CC_{\Bun_G^X}\in \Dmod(\Bun_G)
\end{equation} where $\CC_{\Bun_G^X}$ is the constant $D$-module on $\Bun_G^X$ and $(\pi_{X,\dR})_!$ is the de Rham pushforward.

\subsubsection{Spectral side}
Let $\Gc=\Hc'\times\Hc$ be the Langlands dual group of $G$. Let $\Mc\in \Rep(\Gc)$ be the tensor product of the standard representations of $\Hc$ and $\Hc'$, and let $\Mc^*\in\Rep(\Gc)$ be the dual representation, which is isomorphic to $\Mc$. Let \[\Loc_{\Gc}:=\Map(C_{\dR},\mathbb{B}\Gc)\] be the moduli stack of de Rham $\Gc$-local systems on $C$. 

Let $\Mc_{\univ}^*\in \QCoh(\Loc_{\Gc})\otimes\Dmod(C)$ be the universal $\Gc$-local system on $C$ associated to the representation $\Mc^*$. We use a normalization in which $\Mc_{\univ}^*$ lies in cohomological degree $-2$. Its de Rham cohomology defines \[\Gamma_{\dR}(C,\Mc_{\univ}^*)\in \QCoh(\Loc_{\Gc}).\] This is a derived vector bundle on $\Loc_{\Gc}$ of (virtual) rank $-(2g-2)\dim(\Mc)$. Since $\Mc^*$ is a symplectic self-dual representation of $\Gc$, the de Rham cohomology $\Gamma_{\dR}(C,\Mc_{\univ}^*)[-1]$ carries a symmetric bilinear form. Let $\Loc_{\Gc}^{\circ}\sub \Loc_{\Gc}$ be the open substack where $\Gamma_{\dR}(C,\Mc_{\univ}^*)[-1]$ is concentrated in degree $0$. Following \cite[\S12.8]{BZSV}, one considers the Clifford algebra \[\mathbb{O}_{\Mc}^{\circ}:=\Cl(\Gamma_{\dR}(C,\Mc_{\univ}^*)[-1])\in \operatorname{Alg}(\QCoh(\Loc_{\Gc}^{\circ})).\]

\subsubsection{Comparing the two sides}
The series of papers \cite{gaitsgory2024proof,arinkin2024proof,campbell2024proof,arinkin2024proof2,gaitsgory2024proof2} establishes an equivalence of categories
\[\BL_{G}:\Dmod(\Bun_G) \isom \IndCoh_{\Nilp}(\Loc_{\Gc}).\]
The following conjecture is \cite[Conjecture\,12.8.1]{BZSV} specialized to our setting. It is a geometric Langlands analogue of the Ichino--Ikeda conjecture \cite{ichino2010periods}, which refines the global Gross--Prasad conjecture \cite{Gross_Prasad_1992}:
\begin{conj}\label{conj:main}
The object $\LL_G(\cP_X)|_{\Loc_{\Gc}^{\circ}}$ is an irreducible $\mathbb{O}_{\Mc}^{\circ}$-module in $\IndCoh_{\Nilp}(\Loc_{\Gc}^{\circ})$, that is, there is a canonical isomorphism \[\mathbb{O}_{\Mc}^{\circ}\isom \underline{\End}_{\QCoh(\Loc_{\Gc}^{\circ})}(\LL_G(\cP_X)|_{\Loc_{\Gc}^{\circ}})\in \operatorname{Alg}(\QCoh(\Loc_{\Gc}^{\circ})).\]
\end{conj}

\subsection{Main result}

Let \[j_{\irr}:\Loc_{\Gc}^{\irr}\incl \Loc_{\Gc}\] be the locus of irreducible $\Gc$-local systems. One has $\Loc_{\Gc}^{\irr} \sub \Loc_{\Gc}^{\circ}$. Note that $\pi_0(\Loc_{\Gc}^{\irr})=\pi_0(\Loc_{\Gc})=\pi_1(\Gc)\cong \ZZ/2\ZZ$. Let \[\Loc_{\Gc}^{\irr,0}\sub \Loc_{\Gc}^{\irr}\] be the neutral component. The main result of this paper is the following:
\begin{thm}\label{thm:mainmodule}
  Conjecture \ref{conj:main} holds when restricted to $\Loc_{\Gc}^{\irr,0}$.
\end{thm}

\begin{remark}
Our proof involves computations with the sheaf of twisted differential operators on $\Bun_G$ at the critical level, which corresponds to the pushforward of the structure sheaf along the map $\Op_{\Gc}\to \Loc_{\Gc}$ from the moduli stack of $\Gc$-opers. The image of this map lies in the neutral component $\Loc_{\Gc}^{\irr,0}$, which is why we only prove Conjecture \ref{conj:main} over $\Loc_{\Gc}^{\irr,0}$. However, the remaining component of $\Loc_{\Gc}^{\irr}$ can be treated similarly by considering Kac--Moody localizations of modules other than the vacuum module, corresponding to monodromy-free opers other than regular opers. We plan to return to this in a future version of this article.
\end{remark}

\begin{remark}
The restricted-variation analogue of Theorem \ref{thm:mainmodule} in positive characteristic should follow from the techniques of \cite{gaitsgory2025geometriclanglandspositivecharacteristic}. Combined with the Frobenius trace formalism of \cite{arinkin2022automorphicfunctionstracefrobenius}, this should imply a direct function-field analogue of the Ichino--Ikeda conjecture in \cite{ichino2010periods}, which, as far as we know, remains largely open over both number fields and function fields, except for small values of $n$. Moreover, combining this result with those of \cite{liu2025higherperiodintegralsderivatives,wang2025specialcycleshtukascategorical, wang2026diagonalcyclesshtukasadjoint}, one should be able to deduce a higher-derivative function-field version of the arithmetic Gan--Gross--Prasad conjecture \cite{gan2012symplectic} in the unramified setting, in the spirit of \cite{YZ1}. We plan to address these questions in future work.
\end{remark}

\begin{remark}
For $n<4$, Conjecture \ref{conj:main} follows readily from the results of \cite{feng2025geometric}. For $g<2$, the irreducible locus $\Loc_{\Gc}^{\irr}$ is empty, so the statement of Theorem \ref{thm:mainmodule} is vacuous. We therefore assume $g\geq 2$ and $n\geq 4$ throughout this paper.
\end{remark}

\subsection{Strategy of proof}

The construction in \cite{liu2025higherperiodintegralsderivatives} equips $\LL_G(\cP_X)|_{\Loc_{\Gc}^{\irr,0}}$ with a natural action of the Clifford algebra $\mathbb{O}_{\Mc}^{\circ}|_{\Loc_{\Gc}^{\irr,0}}$ (see Proposition \ref{prop:clifaction}). It therefore suffices to prove the following:

\begin{thm}\label{thm:main}
The object \[\LL_G(\cP_X)|_{\Loc_{\Gc}^{\irr, 0}}[\delta_H-\delta_{N_{\Omega_C^{1/2}}}]\in \QCoh(\Loc_{\Gc}^{\irr, 0})\] is a locally free sheaf in $\QCoh(\Loc_{\Gc}^{\irr, 0})^{\heartsuit}$, with rank equal to $2^{(g-1)\dim(\Mc)}$. Here \begin{itemize} \item $\delta_H=\dim \Bun_H$. \item $\delta_{N_{\Omega_C^{1/2}}}=\dim \Bun_{N_{\Omega_C^{1/2}}}$ where $\Bun_{N_{\Omega_C^{1/2}}}:=\Map(C,\mathbb{B}(N\rtimes\Gm))\times_{\Map(C,\BB\Gm)}\{\Omega_C^{1/2}\}$. Here $\Gm$ acts on $N$ by conjugation through the cocharacter $2\rho\in X_*(T)$.\end{itemize}
\end{thm}

We prove the following two propositions, which together imply Theorem \ref{thm:main}:
\begin{prop}[Proposition \ref{prop:localconstancy}]\label{prop:localconstancyintro}
The sheaf $\LL_G(\cP_X)|_{\Loc_{\Gc}^{\irr}}\in \QCoh(\Loc_{\Gc}^{\irr})$ is cohomologically bounded and has locally free cohomology sheaves.
\end{prop}

\begin{prop}[Proposition \ref{prop:onefiber}]\label{prop:onefiberintro}
There exists $\s\in \Loc_{\Gc}^{\irr,0}(\CC)$ such that the fiber $i_{\s}^*(\LL_G(\cP_X))[\delta_H-\delta_{N_{\Omega_C^{1/2}}}]\in \Vect$ lies in $\Vect^{\heartsuit}$ and has dimension $2^{(g-1)\dim(\Mc)}$, where $i_{\s}:\{\s\}\to \Loc_{\Gc}$ is the corresponding map.
\end{prop}

We will prove Proposition \ref{prop:localconstancyintro} in \S\ref{sec:localconstancy} by showing that the Atiyah class of $\LL_G(\cP_X)$ vanishes locally on $\Loc_{\Gc}^{\irr}$ (Proposition \ref{prop:at=0}), as a consequence of the Plancherel algebra computation in \cite{Braverman_2022} (Corollary \ref{cor:momentglobint}). This implies that $\LL_G(\cP_X)$ admits (not necessarily flat) connections locally on $\Loc_{\Gc}^{\irr}$. We view this result as analogous to part of the multiplicity one argument in \cite[Theorem\,3.1.8]{arinkin2024proof2} (see Remark \ref{rmk:comparemultone} for more details). Our use of the Atiyah class refines the use of singular support theory for ind-coherent sheaves in geometric Langlands. See Remark \ref{rmk:at=0} for a more detailed discussion.

We will prove Proposition \ref{prop:onefiberintro} in \S\ref{sec:period} by computing the period integral of Hecke eigensheaves arising from the quantization of the Hitchin system in \cite{beilinson1991quantization}. For this computation, we use a naturally defined filtration and reduce the problem to a classical limit calculation in \S\ref{sec:classicallimit}. This is similar to the strategy used in the proof of the index formula in \cite[Theorem\,6.1.2.1]{faergeman2025non} (see Remark \ref{rmk:comparefr} for a more detailed comparison).

\subsection{Notation}
In this article, we fix a smooth projective connected curve $C$ over the complex numbers with genus $g\geq 2$. We will fix $G, H, H', \Gc$ as in \S\ref{sec:intro} and use $\frg,\frh,\frh', \frgc$ to denote their Lie algebras. It should be clear from the context that some of our results are valid for more general groups.

By default, we work with presentable $\infty$-categories. Given such a category $\sC$ and objects $x,y$, we write $\Hom(x,y)=\Map_{\sC}(x,y)$ and $\Hom^0(x,y):=\pi_0(\Map_{\sC}(x,y))$.

We use $\Bun_G$ to denote the moduli stack of $G$-bundles. For each finite-dimensional representation $V$ of $G$, we denote by $V_{\univ}\in \QCoh(\Bun_G\times C)$ the universal vector bundle associated to $V$. We use \[\cL_{G, V}:= \det (\Gamma(C,V_{\univ}))\otimes \det(\Gamma(C,\cO_C^{\oplus\dim(V)}))^{-1}\in \Pic(\Bun_G)\] to denote the line bundle on $\Bun_G$ associated to $V$.

We use $\Loc_{\Gc}$ to denote the moduli stack of de Rham $\Gc$-local systems on $C$. Similarly, for each finite-dimensional representation $V$ of $\Gc$, we denote by $V_{\univ}\in \QCoh(\Loc_{\Gc}\times C_{\dR})$ the universal local system associated to $V$, which lies in cohomological degree $-2$.

By default, we work with (derived) stacks over $\CC$ for the \'etale topology. For a smooth algebraic stack $\sY$, we denote by $\sY_{\dR}$ its de Rham stack. We denote by $\Dmod(\sY):=\QCoh(\sY_{\dR})\cong \IndCoh(\sY_{\dR})$ its category of $D$-modules. We have the adjoint pairs \[\ind^l:\QCoh(\sY) \rightleftarrows \Dmod(\sY):\oblv^l\] and \[\ind^r:\IndCoh(\sY)\rightleftarrows \Dmod(\sY): \oblv^r\] such that $\ind^l\cong \ind^r\circ (-\otimes\omega_{\sY})$ where $\omega_{\sY}\in \IndCoh(\sY)$ is the dualizing complex. We use ${\rm D}_{\sY}:=\ind^l(\cO_{\sY})\cong \ind^r(\omega_{\sY})\in \Dmod(Y)$ to denote the sheaf of differential operators on $\sY$.

Given a line bundle $\cL\in \Pic(\sY)$, we also have the category of twisted $D$-modules $\Dmod_{\cL}(\sY)$ which is canonically equivalent to $\Dmod(\sY)$. We have similar adjoint pairs \[\ind^l_{\cL}:\QCoh(\sY)\rightleftarrows \Dmod_{\cL}(\sY):\oblv^l_{\cL}\] and \[\ind^r_{\cL}:\IndCoh(\sY)\rightleftarrows \Dmod_{\cL}(\sY):\oblv^r_{\cL}.\] We define ${\rm D}_{\sY,\cL}:=\ind^l_{\cL}(\cO_{\sY})\in \Dmod_{\cL}(\sY)$. Under the identification $\Dmod_{\cL}(\sY)\cong \Dmod(\sY)$, we have $\ind^l_{\cL}\cong \ind^l\circ(-\otimes \cL^{-1})$ and other similar identifications.

We will fix a square root $\Omega_C^{1/2}$ of the canonical line bundle $\Omega_C$ on $C$, which gives a square root $\cL_{G,\frg}^{1/2}$ of the line bundle $\cL_{G,\frg}$ on $\Bun_G$.  We use $\kappa$ to denote the critical level of $G$, and we use \[\Dmod_{\kappa}(\Bun_G):=\Dmod_{\cL_{G,\frg}^{1/2}}(\Bun_G)\] to denote the category of twisted $D$-modules on $\Bun_G$ at the critical level, and we will sometimes write $\kappa$ instead of $\cL_{G,\frg}^{1/2}$ in our notation.

Given a morphism $f:\sX\to \sY$ of stacks, we write $f_*, f^*, f^!$ to denote the usual functors on $\QCoh$ and $\IndCoh$, and we write $f_{\dR,*}, f_{\dR}^!$ to denote the usual functors on $\Dmod$, whenever they are defined. When $\sY$ is a point, we write $\Gamma(\sX,-):=f_*(-)$ and $\Gamma_{\dR}(\sX,-):=f_{\dR,*}(-)$. We use $\otimes$ (resp. $\otimes^!$) to denote the monoidal operation on $\QCoh(\sY)$ (resp. $\IndCoh(\sY)$). We also use $\otimes$ to denote the action of $\QCoh(\sY)$ on $\IndCoh(\sY)$. We use $\underline{\End}(-)$ to denote the internal endomorphism object in a symmetric monoidal category.

\subsection*{Acknowledgements}
The authors are grateful to their advisor Zhiwei Yun for suggesting the problem and for constant encouragement and support.

\subsection*{AI usage}
The authors originally had the idea of using similar methods to compute the Euler characteristics of the stalks of $\LL_G(\cP_X)$. Wang then prompted GPT-6\,Pro with this idea and asked whether stronger statements beyond Euler characteristics could be established. This led to Proposition \ref{prop:localconstancyintro}. We also used GPT-6 and Claude\,Opus\,5.5 for literature searches, most notably for Proposition \ref{prop:HNstrat}, and for spotting typos and polishing the grammar. The article was written by humans, and the authors take full responsibility for any errors.

\section{Local constancy}\label{sec:localconstancy}

\subsection{Main result}
In this section, we prove the following result:

\begin{prop}[Proposition \ref{prop:localconstancyintro}]\label{prop:localconstancy}
The sheaf $\LL_G(\cP_X)|_{\Loc_{\Gc}^{\irr}}\in \QCoh(\Loc_{\Gc}^{\irr})$ is cohomologically bounded and has locally free cohomology sheaves.
\end{prop}

We will deduce Proposition \ref{prop:localconstancy} from the following two lemmas:

\begin{lemma}\label{lem:coherentcoh}
The object $\LL_G(\cP_X)|_{\Loc_{\Gc}^{\irr}}$ is a perfect complex.
\end{lemma}

\begin{lemma}\label{lem:localconnection}
The object $\LL_G(\cP_X)|_{\Loc_{\Gc}^{\irr}}$ admits (not necessarily flat) connections locally on $\Loc_{\Gc}^{\irr}$.
\end{lemma}

\begin{remark}\label{rmk:comparemultone}
Lemma \ref{lem:localconnection} can be compared with \cite[Theorem\,3.1.8]{arinkin2024proof2}, which equips $\LL_G(\cP_{\mathrm{Whit}})|_{\Loc_{\Gc}^{\irr}}$ with a flat connection, where $\cP_{\mathrm{Whit}}\in\Dmod(\Bun_G)$ is the Whittaker period sheaf.\footnote{It is denoted by $\mathrm{Poinc}_{G,!}^{\mathrm{Vac},\mathrm{glob}}$ in \cite{gaitsgory2024proof}.} Our setting differs from that of \textit{loc. cit.} in two respects. First, we only obtain local, non-canonical, and not necessarily flat connections on $\Loc_{\Gc}^{\irr}$, whereas the connection constructed in \textit{loc. cit.} is canonical and flat. Second, the proof strategies are quite different: we work directly with structures on the period sheaf, whereas \textit{loc. cit.} uses ambidexterity to reduce to the study of the period associated with opers. Despite these differences, we view both results as establishing certain local constancy properties needed for multiplicity one arguments. Our method can also be applied to equip $\LL_G(\cP_{\mathrm{Whit}})|_{\Loc_{\Gc}^{\irr}}$ with local, non-canonical, and not necessarily flat connections.

\end{remark}

It is easy to see that Lemma \ref{lem:coherentcoh} and Lemma \ref{lem:localconnection} imply Proposition \ref{prop:localconstancy}. Lemma \ref{lem:localconnection} follows immediately from Proposition \ref{prop:at=0}, which we prove below. We first prove Lemma \ref{lem:coherentcoh}.

\begin{proof}[Proof of Lemma \ref{lem:coherentcoh}]
Let \begin{equation}\label{eq:buncusp} j_{\cusp}:\Bun_G^{\cusp}\to \Bun_G\end{equation} be the open immersion of the open substack where the Harder--Narasimhan filtration has slope $\lambda\in X_*(T)^+_{\QQ}$ with $\lambda \leq (2g-2)\rho$ where $\rho\in X_*(T)^+_{\QQ}$ is the half-sum of positive coroots as defined in \cite[\S B.1.2]{drinfeld2016geometric}. Here the order is defined such that $\lambda_1\leq \lambda_2$ if $\lambda_2-\lambda_1$ is a linear combination of simple coroots of $G$ with non-negative coefficients. By \cite[Proposition\,1.4.6]{drinfeld2016geometric}, any object of $\Dmod(\Bun_G)_{\cusp}$ is a clean extension of its restriction to $\Bun_G^{\cusp}$. Denote $\cF := j_{\cusp}^*\cP_X\in \Dmod(\Bun_G^{\cusp})$. Let
\[\LL_{G}^{\coarse}:\Dmod(\Bun_G)\to \QCoh(\Loc_{\Gc})\]
denote the coarse Langlands functor considered in \cite[\S1.4]{gaitsgory2024proof}, defined as the composition of $\LL_G$ with the natural functor $\IndCoh_{\Nilp}(\Loc_{\Gc})\to \QCoh(\Loc_{\Gc})$. We then have
\[\LL_G(\cP_X)|_{\Loc_{\Gc}^{\irr}}\cong \LL_G^{\coarse}((j_{\cusp})_!\cF)|_{\Loc_{\Gc}^{\irr}}.\]

First, we show that $\LL_G^{\coarse}((j_{\cusp})_!\cF)$ has bounded cohomological degrees. Indeed, by \cite[Proposition\,1.8.2, Theorem\,2.1.2]{gaitsgory2024proof}, we know that $\LL_G^{\coarse}$ has bounded cohomological amplitude. Therefore, it suffices to show that $(j_{\cusp})_!\cF$ has bounded cohomological degrees. This follows from the fact that $\cF\in \Dmod(\Bun_G^{\cusp})$ has bounded cohomological degrees and that $(j_{\cusp})_!$ has bounded cohomological amplitude, where the latter follows from the argument in \cite[\S2.2.4]{gaitsgory2024proof}: it is clear that $(j_{\cusp})_!$ is right $t$-exact, and the left cohomological amplitude is bounded by the minimal number of affine charts covering a smooth cover of $\Bun_G^{\cusp}$.

Assume $\cF\in \Dmod(\Bun_G^{\cusp})^{\leq N}$ for some $N\in \ZZ$. Write $\Bun_G^{\cusp}=[Z/K]$ where $Z$ is a smooth scheme and $K$ is a smooth affine algebraic group, and let $f:Z\to \Bun_G^{\cusp}$ be the quotient map. We have
\[\cF\cong \colim(\cdots f_!f^!f_!f^!\cF\rightrightarrows f_!f^!\cF).\]
Let $\cF_d$ denote the geometric realization of the $d$-skeleton of this simplicial object. Then $\operatorname{Cone}(\cF_d\to \cF)\in \Dmod(\Bun_G^{\cusp})^{\leq N-d}$. Since all $\cF_d$ are compact, we know that $\LL_G^{\coarse}((j_{\cusp})_!\cF_d)|_{\Loc_{\Gc}^{\irr}}$ are perfect complexes. Combining this with the cohomological boundedness of $\LL_G^{\coarse}\circ (j_{\cusp})_!$, we obtain that $\LL_G^{\coarse}((j_{\cusp})_!\cF)|_{\Loc_{\Gc}^{\irr}}$ has coherent cohomology sheaves and hence is a perfect complex.

\end{proof}

\subsection{Atiyah class}

\subsubsection{First infinitesimal neighborhood}\label{sec:firstinf}
We first recall the construction of the first infinitesimal neighborhood of a morphism in \cite[Chapter\,9,\,\S5.1]{gaitsgory2020study}. Given a morphism between prestacks \[f:{\scr Y}\to {\scr Z}\] where $\sY$ and $\sZ^{\wedge}_{\sY}:=\sZ\times_{\sZ_{\dR}}\sY_{\dR}$ are locally almost of finite type (laft) and admit deformation theory in the sense of \cite[Chapter\,1,\,\S7.1]{gaitsgory2020study}, consider the associated map on tangent complexes \footnote{$T_{\sY}\in \IndCoh(\sY)$ should be regarded as the \emph{Serre} dual of the cotangent complex $\LL_{\sY}$, rather than the linear dual. See \cite[Chapter\,1,\,\S3.5.7]{gaitsgory2020study}.} \[T_{\sY/\sZ}\to T_{\sY}.\] The construction in \cite[Chapter\,8,\,\S5.1]{gaitsgory2020study} produces a square-zero extension\footnote{We will sometimes write $\sY^{(1)}$ rather than $\sY^{(1)}_{\sZ}$ when $\sZ$ is clear from the context.} \[f^{(1)}:\sY\to \sY^{(1)}_{\sZ}\] equipped with a map \[\sY^{(1)}_{\sZ}\to \sZ.\] This is the first infinitesimal neighborhood of the morphism $f$. Moreover, \cite[Chapter\,8,\,Proposition\,6.4.2]{gaitsgory2020study} gives a fiber sequence in $\IndCoh(\sY^{(1)}_{\sZ})$:\footnote{The functor $f^{(1)}_*$ is defined since $f^{(1)}$ is inf-schematic, and we use the $\IndCoh$-sheaf theory in \cite[Chapter\,3,\,Theorem\,5.4.3]{gaitsgory2020study}. The same applies to the other $*$-pushforward functors in this section.}
\begin{equation}\label{eq:infseq}
(f^{(1)})_* T_{\sY/\sZ}\to (f^{(1)})_* \omega_{\sY}\to \omega_{\sY^{(1)}_{\sZ}}.
\end{equation}

This construction is functorial in $f$: Given a sequence of morphisms \[\sY\xrightarrow{f} \sZ_1\xrightarrow{g}\sZ_2,\] there is an induced map \[i_g:\sY^{(1)}_{\sZ_1}\to \sY^{(1)}_{\sZ_2}\] and a commutative diagram in $\IndCoh(\sY^{(1)}_{\sZ_2})$:
\begin{equation}\label{eq:firstinffunct}\begin{tikzcd}(i_g\circ f^{(1)}_{\sZ_1})_* T_{\sY/\sZ_1} \ar[r] \ar[d] & (i_g\circ f^{(1)}_{\sZ_1})_* \omega_{\sY} \ar[r] \ar[d] & (i_g)_* \omega_{\sY^{(1)}_{\sZ_1}} \ar[d] \\(f^{(1)}_{\sZ_2})_* T_{\sY/\sZ_2} \ar[r] & (f^{(1)}_{\sZ_2})_* \omega_{\sY} \ar[r] & \omega_{\sY^{(1)}_{\sZ_2}}  \end{tikzcd}.\end{equation}

The construction also commutes with base change by \cite[Chapter\,8,\,Proposition\,5.4.3]{gaitsgory2020study}: Given a Cartesian diagram 
\[\begin{tikzcd}
\widetilde{\sY} \ar[r, "\widetilde{f}"] \ar[d, "g_{\sY}"] & \widetilde{\sZ} \ar[d, "g_{\sZ}"] \\ \sY \ar[r, "f"] & \sZ
\end{tikzcd}\] in which $\sY,\widetilde{\sY},\sZ^{\wedge}_{\sY},\widetilde{\sZ}^{\wedge}_{\widetilde{\sY}}$ are laft and admit deformation theory, one has a commutative diagram with Cartesian squares:
\[\begin{tikzcd}\widetilde{\sY}\ar[r, "\widetilde{f}^{(1)}"]\ar[d, "g_{\sY}"] & \widetilde{\sY}^{(1)}_{\widetilde{\sZ}} \ar[r] \ar[d, "g_{\sY}^{(1)}"] & \widetilde{\sZ} \ar[d, "g_{\sZ}"] \\ \sY \ar[r, "f^{(1)}"] & \sY^{(1)}_{\sZ} \ar[r] & \sZ\end{tikzcd}\] and a commutative diagram in $\IndCoh(\widetilde{\sY}^{(1)}_{\widetilde{\sZ}})$:
\begin{equation}\label{eq:firstinfbc}\begin{tikzcd}
( \widetilde{f}^{(1)})_* T_{\widetilde{\sY}/\widetilde{\sZ}} \ar[r] \ar[d, "\sim"] & (\widetilde{f}^{(1)})_* \omega_{\widetilde{\sY}} \ar[r] \ar[d, "\sim"] &  \omega_{\widetilde{\sY}^{(1)}_{\widetilde{\sZ}}} \ar[d, "\sim"] \\ (g_{\sY}^{(1)})^!(f^{(1)})_* T_{\sY/\sZ} \ar[r] & (g_{\sY}^{(1)})^!(f^{(1)})_* \omega_{\sY} \ar[r] & (g_{\sY}^{(1)})^!\omega_{\sY^{(1)}_{\sZ}} \end{tikzcd}.\end{equation}

\subsubsection{General definition}\label{sec:atgeneral}
Applying the construction of \S\ref{sec:firstinf} to the diagonal map \[\Delta_{\scr Y}:{\scr Y}\to {\scr Y}\times {\scr Y},\] one obtains the first infinitesimal neighborhood of the diagonal:
  \[\sY\xrightarrow{\D_{\sY}^{(1)}} \sY^{(1)}\xrightarrow{(p_1^{(1)},p_2^{(1)})} \sY\times\sY\] and a fiber sequence in $\IndCoh(\sY^{(1)})$:
\begin{equation}\label{eq:atiyahclassker}
(\D_{\sY}^{(1)})_* T_{\sY}[-1]\to (\D_{\sY}^{(1)})_* \omega_{\sY}\to \omega_{\sY^{(1)}}.
\end{equation}
Regarding these as kernels gives a natural transformation between endofunctors on $\IndCoh(\sY)$:\footnote{Since the maps $\sY^{(1)}\to (\sY\times\sY)^{\wedge}_{\sY}\rightrightarrows \sY$ are inf-schematic, the kernel perspective is always valid.}
\begin{equation}\label{eq:atiyahclass}
\At_{\sY}: T_{\sY}[-1]\otimes^!-\to \id.
\end{equation}
We call \eqref{eq:atiyahclass} (or \eqref{eq:atiyahclassker}) the \emph{Atiyah class} of $\sY$. For $\cF\in \IndCoh(\sY)$, we denote by \[\At_{\sY}(\cF): T_{\sY}[-1]\otimes^!\cF \to \cF\] the induced morphism.

For $\cF\in \IndCoh(\sY)$, a \emph{connection} on $\cF$ is a null-homotopy of $\At_{\sY}(\cF)$. To see that this is compatible with the usual definition of a connection, one notes that \eqref{eq:atiyahclassker} gives a fiber sequence in $\IndCoh(\sY)$:
\[T_{\sY}[-1]\otimes^!\cF \to \cF\to (p_2^{(1)})_*(p_1^{(1)})^! \cF .\] Therefore, a null-homotopy of $\At_{\sY}(\cF)$ is equivalent to a splitting of the last map, which by adjunction is equivalent to a map $(p_1^{(1)})^!\cF\to (p_2^{(1)})^!\cF$ in $\IndCoh(\sY^{(1)})$ whose restriction to $\sY$ is the identity. This is exactly the usual definition of a connection.

\subsubsection{Atiyah class of $\Loc_{\Gc}$}
Applying the construction in \S\ref{sec:atgeneral} to $\sY=\Loc_{\Gc}$, we obtain the Atiyah class of $\Loc_{\Gc}$ as a natural transformation between endofunctors on $\IndCoh(\Loc_{\Gc})$:
\begin{equation}\label{eq:at}
\At: \Gamma_{\dR}(C,\frgc_{\univ})[-2]\otimes-\cong  T_{\Loc_{\Gc}}[-1]\otimes^!-\to \id,
\end{equation}
where we regard $\Gamma_{\dR}(C,\frgc_{\univ})\in \QCoh(\Loc_{\Gc})$. It satisfies $\Gamma_{\dR}(C,\frgc_{\univ})[-2]\otimes\omega_{\Loc_{\Gc}}\cong T_{\Loc_{\Gc}}[-1]$.

We also need the variant where one takes $\sY = \Loc_{\Gc}\times C_{\dR}$, which gives an Atiyah class as a natural transformation between endofunctors on $\IndCoh(\Loc_{\Gc}\times C_{\dR})$:\footnote{Note that under the identification $\Dmod(\sY)\cong\QCoh(\sY_{\dR})$, the usual tensor product $\otimes$ on $\QCoh(\sY_{\dR})$ corresponds to the $!$-tensor product on $\Dmod(\sY)$, and we use this tensor product by default.}
\begin{equation}\label{eq:atC}
\At_C: (\Gamma_{\dR}(C,\frgc_{\univ})[-2]\boxtimes\omega_{C_{\dR}})\otimes-\to \id,
\end{equation} where we are using the fact that $T_{\Loc_{\Gc}\times C_{\dR}}\cong T_{\Loc_{\Gc}}\boxtimes \omega_{C_{\dR}}$. It is straightforward to verify that $\At_C\cong \At\otimes\id$ as natural transformations of endofunctors of $\IndCoh(\Loc_{\Gc}\times C_{\dR})\cong\IndCoh(\Loc_{\Gc})\otimes\Dmod(C)$.

\subsection{Local Atiyah class}


\subsubsection{Spherical Hecke category}
For an affine scheme $X$, we write $L^+X$ and $LX$ for its jet space and loop space, respectively, and $(L^+X)_C$ and $(LX)_C$ for their respective versions over $C_{\dR}$.

Define the spherical Hecke category \[\Sph_G:=\Dmod(L^+G\bs LG/L^+G)^{\ren}\cong \IndCoh(\mathbb{B}\Gc\times_{[\frgc/\Gc]} \mathbb{B}\Gc).\] Here the superscript $\ren$ stands for the renormalization as in \cite[\S1.5.2]{arinkin2024proof}. This construction extends to families over $C$ and defines a monoidal category over $\Dmod(C)$ \[\Sph_{G,C}:=\Dmod((L^+G)_C\bs (LG)_C/(L^+G)_C)^{\ren}\cong \IndCoh(\Loc_{\Gc,C}^{\reg}\times_{\Loc_{\Gc,C}^{\mer}}\Loc_{\Gc,C}^{\reg})\] where $\Loc_{\Gc,C}^{\reg}$ (and $\Loc_{\Gc,C}^{\mer}$) are informally defined as the moduli stacks classifying a point $c\in C_{\dR}$ together with a $\Gc$-local system on the formal disk (and punctured disk) around $c$. We refer to \cite[\S B.7]{arinkin2024proof} for a rigorous definition.\footnote{In \textit{loc. cit.}, these stacks are defined over $\Ran$, while our $\Loc_{\Gc,C}^{\reg}$ and $\Loc_{\Gc,C}^{\mer}$ are their restrictions to the diagonal $C\sub \Ran$.} We use $\Sat:\Rep(\Gc)\to \Sph_G$ and $\Sat_{C}:\Rep(\Gc)\to \Sph_{G,C}$ to denote the obvious functors.

Moreover, one also considers the spherical Hecke category with loop rotation \begin{equation}\label{eq:dersat}\Sph_{G,\hbar}:=\Dmod((L^+G\rtimes\Aut(D))\bs (LG\rtimes \Aut(D))/ (L^+G\rtimes\Aut(D)))^{\ren}\cong \Mod_{U_{\hbar}(\frgc^{\shear})}^{\Gc}\end{equation} where $\Aut(D)$ is the group scheme of automorphisms of the formal disk $D=\Spec(\CC[\![t]\!])$ fixing the origin, and $U_{\hbar}(\frgc^{\shear})$ is the Rees version of the universal enveloping algebra of $\frgc$. Here, the shearing $\shear$ indicates that $\frgc$ is placed in cohomological degree $2$. This is \cite[Theorem\,5]{BF}. We use $\Sat_{\hbar}:\Rep(\Gc)\to \Sph_{G,\hbar}$ to denote the natural functor. 

We use $*$ to denote the monoidal operation given by convolution on $\Sph_G$, $\Sph_{G,C}$ and $\Sph_{G,\hbar}$, and we use the same notation $*$ to denote the action of $\Sph_{G,C}$ on $\Dmod(\Bun_G)$,\footnote{By the ``action'' of $\Sph_{G,C}$ on $\Dmod(\Bun_G)$, we mean the composition of the usual Hecke action with integration over $C$. Note that this ``action'' is not compatible with the monoidal structure on $\Sph_{G,C}$ over $\Dmod(C)$.} $\Dmod(\Bun_G\times C)$ (as well as the equivalent action on the spectral side) understood in the usual way. Note that the compatibility of the actions under the equivalence $\LL_G$ is a special case of \cite[Proposition\,18.3.8]{arinkin2024proof} which proves the compatibility over the entire Ran space.

\subsubsection{Local Atiyah classes}

Koszul duality gives an equivalence of categories \[\Sph_G\cong \QCoh([\frgc^*/\Gc])^{\shear}\] where the shearing $\shear$ places $\frgc^*$ in cohomological degree $-2$. In particular, this description gives a tautological morphism \begin{equation}\label{eq:atloc}\At^{\loc}: \Sat(\frgc)\to \Sat(\CC)[2].\end{equation} We call this the \emph{local Atiyah class}.

We have a unique lift of the local Atiyah class to a morphism \begin{equation}\label{eq:atlochbar}\At^{\loc}_{\hbar}: \Sat_{\hbar}(\frgc)\to \Sat_{\hbar}(\CC)[2]\end{equation} also given tautologically by \eqref{eq:dersat}.

The following lemma gives a unique lift of the local Atiyah class to a family over the curve:
\begin{lemma}
\label{lem:atlocC}
There is a contractible space of choices of a morphism \begin{equation}\label{eq:atlocC}\At^{\loc}_C: \Sat_{C}(\frgc)\to \Sat_{C}(\CC)[2]\end{equation} whose fiber at any point $c\in C$ is the local Atiyah class $\At^{\loc}$.
\end{lemma}
\begin{proof}
The existence follows from pulling back \eqref{eq:atlochbar} on the automorphic side. Again on the automorphic side, one has \[\Hom^{*}(\Sat_{C}(\frgc), \Sat_{C}(\CC)[2])\cong \Hom^{*}(\Sat_{\hbar}(\frgc), \Sat_{\hbar}(\CC)[2])\otimes_{\CC[\hbar]}H^*_{\dR}(C,\CC).\] Since $\Hom^{i}(\Sat_{\hbar}(\frgc), \Sat_{\hbar}(\CC)[2])=0$ for $i<0$, we know that \[\Hom^{0}(\Sat_{C}(\frgc), \Sat_{C}(\CC)[2])\cong \Hom^0(\Sat_{\hbar}(\frgc),\Sat_{\hbar}(\CC)[2])\cong \Hom^0(\Sat(\frgc),\Sat(\CC)[2]).\] This implies the uniqueness of \eqref{eq:atlocC}.
\end{proof}


By abuse of notation, we also write $\At^{\loc}_C$ for the induced natural transformation between endofunctors on $\IndCoh(\Loc_{\Gc}\times C_{\dR})$:
\begin{equation}\label{eq:atlocCend}
\At^{\loc}_C: \frgc_{\univ}[-2]\otimes -\to \id.
\end{equation}

\subsubsection{An alternative description of the local Atiyah class}
We first give an alternative description of the local Atiyah class \eqref{eq:atloc}. Applying the construction of \S\ref{sec:firstinf} to the diagonal map \[\D_{\BB\Gc}:\BB\Gc\to \BB\Gc\times_{[\frgc/\Gc]}\BB\Gc,\] we get a square-zero extension \[\BB\Gc\xrightarrow{\D_{\BB\Gc}^{(1)}} (\BB\Gc)^{(1)} \xrightarrow{(p_1^{(1)},p_2^{(1)})} \BB\Gc\times_{[\frgc/\Gc]}\BB\Gc\] and a fiber sequence in $\IndCoh(\BB\Gc^{(1)})$:
\begin{equation}\label{eq:atiyahclasskerloc}
(\D_{\BB\Gc}^{(1)})_* (T_{\BB\Gc/(\BB\Gc\times_{[\frgc/\Gc]}\BB\Gc)})\xrightarrow{\At^{\loc}} (\D_{\BB\Gc}^{(1)})_* \omega_{\BB\Gc}\to \omega_{\BB\Gc^{(1)}}
.\end{equation} The first map is a lift of the local Atiyah class \eqref{eq:atloc} in the sense that after applying $(p_1^{(1)},p_2^{(1)})_*$ to \eqref{eq:atiyahclasskerloc} and using the identification \[(\D_{\BB\Gc})_*T_{\BB\Gc/(\BB\Gc\times_{[\frgc/\Gc]}\BB\Gc)}\cong \Sat(\frgc)[-2],\] it becomes \eqref{eq:atloc}.

We can also apply the construction in \S\ref{sec:firstinf} to the diagonal map \[\D_{\Loc_{\Gc,C}^{\reg}}:\Loc_{\Gc,C}^{\reg}\to \Loc_{\Gc,C}^{\reg}\times_{\Loc_{\Gc,C}^{\mer}}\Loc_{\Gc,C}^{\reg},\] which gives a square-zero extension \[\Loc_{\Gc,C}^{\reg}\xrightarrow{\D_{\Loc_{\Gc,C}^{\reg}}^{(1)}} \Loc_{\Gc,C}^{\reg,(1)} \to \Loc_{\Gc,C}^{\reg}\times_{\Loc_{\Gc,C}^{\mer}}\Loc_{\Gc,C}^{\reg}\] and a fiber sequence in $\IndCoh(\Loc_{\Gc,C}^{\reg,(1)})$:
 \begin{equation}\label{eq:atiyahclasskerlocC}
(\D_{\Loc_{\Gc,C}^{\reg}}^{(1)})_* (T_{\Loc_{\Gc,C}^{\reg}/(\Loc_{\Gc,C}^{\reg}\times_{\Loc_{\Gc,C}^{\mer}}\Loc_{\Gc,C}^{\reg})})\xrightarrow{\At^{\loc}_C} (\D_{\Loc_{\Gc,C}^{\reg}}^{(1)})_* \omega_{\Loc_{\Gc,C}^{\reg}}\to \omega_{(\Loc_{\Gc,C}^{\reg})^{(1)}}.
\end{equation}
where the first map is a lift of the Atiyah class \eqref{eq:atlocC} in the same sense under the isomorphism \[(\D_{\Loc_{\Gc,C}^{\reg}})_*T_{\Loc_{\Gc,C}^{\reg}/(\Loc_{\Gc,C}^{\reg}\times_{\Loc_{\Gc,C}^{\mer}}\Loc_{\Gc,C}^{\reg})}\cong \Sat_{C}(\frgc)[-2].\] Given this isomorphism, the lifting property follows immediately from Lemma \ref{lem:atlocC} and the base change property of first infinitesimal neighborhoods in \eqref{eq:firstinfbc}. To see this isomorphism, we only need to show that there is an isomorphism in $\IndCoh(\Loc_{\Gc,C}^{\reg})\cong \IndCoh(\BB\Gc\times C_{\dR})$:
\[ T_{\Loc_{\Gc,C}^{\reg}/\Loc_{\Gc,C}^{\mer}}\cong \frgc[-1]\boxtimes \omega_{C_{\dR}}.
\]
This follows from \cite[Proposition\,6.1.2]{rozenblyum2021connectionsmodulispacesinfinitesimal} by pulling back along the map $\BB\Gc\to \Loc_{\Gc}$ given by the trivial $\Gc$-local system on $C$.

\subsection{Comparing local and global Atiyah classes}
Now we compare the local Atiyah class \eqref{eq:atlocC} with the global Atiyah class \eqref{eq:atC}. The comparison is given by the following proposition:

\begin{prop}\label{prop:atloc=at}
One has the following commutative diagram of natural transformations between endofunctors on $\IndCoh(\Loc_{\Gc}\times C_{\dR})$:
\[\begin{tikzcd}
  \frgc_{\univ}[-2]\otimes - \ar[r, "\At_C^{\loc}"] \ar[d, "\Gys\otimes-"] & \id \ar[d, "\sim"]\\
  (\Gamma_{\dR}(C,\frgc_{\univ})[-2]\boxtimes \omega_{C_{\dR}})\otimes- \ar[r, "\At_C"] & \id
\end{tikzcd}\] where $\Gys:\frgc_{\univ}[-2]\to \Gamma_{\dR}(C,\frgc_{\univ})[-2]\boxtimes \omega_{C_{\dR}}\in \QCoh(\Loc_{\Gc}\times C_{\dR})$ is the Gysin map induced by the $((p_{C,\dR})_*,p_{C,\dR}^!)$-adjunction, where $p_C:C\to *$ is the morphism to a point.
\end{prop}
\begin{proof}
Since all the endofunctors are linear over $\Dmod(C)$, it suffices to compare the kernels of the natural transformations as objects in $\IndCoh((\Loc_{\Gc}\times \Loc_{\Gc})^{\wedge}_{\Loc_{\Gc}}\times C_{\dR})$. 

We first describe the local Atiyah class $\At_C^{\loc}$ in terms of kernels. Let $\Loc_{\Gc,C}^{\mer,\glob}$ be the moduli stack of a point $c\in C_{\dR}$ and a $\Gc$-local system on $C\setminus c$. Using \cite[Lemma\,E.7.2]{arinkin2024proof}, we have a Cartesian diagram
\begin{equation}\label{eq:loclocglob}\begin{tikzcd}
\Loc_{\Gc}\times C_{\dR} \ar[r, "\D_{\Loc_{\Gc}\times C_{\dR}}^{\loc}"] \ar[d] & (\Loc_{\Gc}\times C_{\dR})\times_{\Loc_{\Gc,C}^{\mer,\glob}}(\Loc_{\Gc}\times C_{\dR}) \ar[d]\\
\Loc_{\Gc,C}^{\reg} \ar[r, "\D_{\Loc_{\Gc,C}^{\reg}}"] & \Loc_{\Gc,C}^{\reg}\times_{\Loc_{\Gc,C}^{\mer}}\Loc_{\Gc,C}^{\reg}
\end{tikzcd}.\end{equation}
Applying the construction of \S\ref{sec:firstinf} to the upper horizontal map, we obtain a square-zero extension \[\Loc_{\Gc}\times C_{\dR}\xrightarrow{(\D_{\Loc_{\Gc}\times C_{\dR}}^{\loc})^{(1)}} (\Loc_{\Gc}\times C_{\dR})^{(1)}_{(\Loc_{\Gc}\times C_{\dR})\times_{\Loc_{\Gc,C}^{\mer,\glob}}(\Loc_{\Gc}\times C_{\dR})}\to (\Loc_{\Gc}\times C_{\dR})\times_{\Loc_{\Gc,C}^{\mer,\glob}}(\Loc_{\Gc}\times C_{\dR})\] and a fiber sequence in $\IndCoh((\Loc_{\Gc}\times C_{\dR})^{(1)})$: 

\[
\adjustbox{max width=\linewidth}{$\displaystyle
(\D_{\Loc_{\Gc}\times C_{\dR}}^{\loc})^{(1)}_*T_{(\Loc_{\Gc}\times C_{\dR})/(\Loc_{\Gc}\times C_{\dR})\times_{\Loc_{\Gc,C}^{\mer,\glob}}(\Loc_{\Gc}\times C_{\dR})}\xrightarrow{\At_C^{\loc}} (\D_{\Loc_{\Gc}\times C_{\dR}}^{\loc})^{(1)}_*\omega_{\Loc_{\Gc}\times C_{\dR}}\to \omega_{(\Loc_{\Gc}\times C_{\dR})^{(1)}_{(\Loc_{\Gc}\times C_{\dR})\times_{\Loc_{\Gc,C}^{\mer,\glob}}(\Loc_{\Gc}\times C_{\dR})}}
$}
\]

where the first map induces the local Atiyah class $\At_C^{\loc}$ at the level of kernels. This follows by unwinding the definition of the Hecke action and the description of the local Atiyah class in \eqref{eq:atiyahclasskerlocC}, together with the base change property of first infinitesimal neighborhoods \eqref{eq:firstinfbc} applied to \eqref{eq:loclocglob}.

Now we compare the local Atiyah class $\At_C^{\loc}$ with the global Atiyah class $\At_C$. Consider the sequence of morphisms \[\Loc_{\Gc}\times C_{\dR}\xrightarrow{\D_{\Loc_{\Gc}\times C_{\dR}}^{\loc}} (\Loc_{\Gc}\times C_{\dR})\times_{\Loc_{\Gc,C}^{\mer,\glob}}(\Loc_{\Gc}\times C_{\dR})\xrightarrow{\D_{\Loc_{\Gc}\times C_{\dR}}^{\loc\to\glob}} \Loc_{\Gc}\times\Loc_{\Gc}\times C_{\dR}\] and apply the functoriality of the first infinitesimal neighborhood. We obtain an induced map \[i_{\D_{\Loc_{\Gc}\times C_{\dR}}^{\loc\to\glob}}:(\Loc_{\Gc}\times C_{\dR})^{(1)}_{(\Loc_{\Gc}\times C_{\dR})\times_{\Loc_{\Gc,C}^{\mer,\glob}}(\Loc_{\Gc}\times C_{\dR})}\to (\Loc_{\Gc}\times C_{\dR})^{(1)}_{\Loc_{\Gc}\times\Loc_{\Gc}\times C_{\dR}}\] and also a commutative diagram as in \eqref{eq:firstinffunct}. This reduces the problem to showing the following commutative diagram in $\IndCoh(\Loc_{\Gc}\times C_{\dR})$:
\[\begin{tikzcd}
T_{(\Loc_{\Gc}\times C_{\dR})/(\Loc_{\Gc}\times C_{\dR})\times_{\Loc_{\Gc,C}^{\mer,\glob}}(\Loc_{\Gc}\times C_{\dR})} \ar[d] \ar[r, "\sim"] & \frgc_{\univ}[-2]\otimes\omega_{\Loc_{\Gc}\times C_{\dR}} \ar[d, "\Gys\otimes \id"] \\
T_{(\Loc_{\Gc}\times C_{\dR})/(\Loc_{\Gc}\times\Loc_{\Gc}\times C_{dR})} \ar[r, "\sim"] & (\Gamma_{\dR}(C,\frgc_{\univ})[-2]\boxtimes\omega_{C_{\dR}})\otimes\omega_{\Loc_{\Gc}\times C_{\dR}} \end{tikzcd}.\] This follows from \cite[Proposition\,5.3.3(b)]{rozenblyum2021connectionsmodulispacesinfinitesimal} after a shifting.

\end{proof}

\subsection{\texorpdfstring{$L$-observables}{L-observables}}

We will prove the following result in this section:
\begin{prop}\label{prop:at=0}
The morphism $\At(\LL_G(\cP_X)):\Gamma_{\dR}(C,\frgc_{\univ})[-2]\otimes \LL_G(\cP_X) \to \LL_G(\cP_X)$ is null-homotopic locally on $\Loc_{\Gc}^{\irr}$.
\end{prop}

This immediately implies Lemma \ref{lem:localconnection}. The proof of Proposition \ref{prop:at=0} will be given in \S\ref{sec:proofat=0} after some preparation.

\begin{remark}\label{rmk:at=0}
    The Atiyah class also plays a role in the theory of singular support for ind-coherent sheaves developed in \cite{arinkin2015singular}, which has become a standard tool in the study of geometric Langlands. In the theory of singular support, one studies the Atiyah class only up to nilpotence. For example, on a quasi-smooth scheme $S$, a sheaf $\cF\in\Coh(S)$ has zero singular support if and only if its Atiyah class $\At(\cF)$ is nilpotent. In contrast, Proposition \ref{prop:at=0} addresses the finer question of where the Atiyah class vanishes identically. To the best of our knowledge, this perspective is new in the context of geometric Langlands.
\end{remark}


\subsubsection{Plancherel algebra}
Consider the relative derived Satake category \[\Sph_{X}^{\mathrm{nonren}}:=\Shv(LX/L^+G)\] where $\Shv(LX/L^+G)$ is the category of ind-holonomic $D$-modules on $LX/L^+G$ as defined in \cite[\S B.7]{BZSV}.\footnote{This category is denoted $\mathrm{SHV}_s^!(LX/L^+G)$ in \cite[\S B.7]{BZSV}.} Let $\Sph_{G}^{\mathrm{nonren}}:=\Shv(L^+G\bs LG/L^+G)\sub \Sph_G$ be the non-renormalized spherical Hecke category. Then $\Sph_{X}^{\mathrm{nonren}}$ is a module category over $\Sph_{G}^{\mathrm{nonren}}$. It has a loop rotation version \[\Sph_{X,\hbar}^{\mathrm{nonren}}:=\Shv(LX/(L^+G\rtimes \Aut(D)))\] which is a module category over $\Sph_{G,\hbar}^{\mathrm{nonren}}$.

Consider the map $i_X:L^+X/L^+G\to LX/L^+G$ and let $\delta_X:=(i_X)_*\omega_{L^+X/L^+G}\in \Sph_{X}^{\mathrm{nonren}}$, and similarly define $\delta_{X,\hbar}\in \Sph_{X,\hbar}^{\mathrm{nonren}}$.  The Plancherel algebra is defined as \[\PL_X:=\underline{\End}_{\Rep(\Gc)}(\delta_X)\in \operatorname{Alg}(\Rep(\Gc))\] and it admits a quantization \[\PL_{X,\hbar}:=\underline{\End}_{\Rep(\Gc)}(\delta_{X,\hbar})\in \operatorname{Alg}(\Rep(\Gc)).\] 

For our choice of $X$, \cite[Theorem\,2.2.1, Theorem\,3.1.1]{Braverman_2022} gives isomorphisms \begin{equation}\label{eq:plX}
\PL_X\cong \cO(\Mc)^{\shear}
\end{equation} and \begin{equation}\label{eq:plXhbar}
\PL_{X,\hbar}\cong \cO(\Mc)^{\shear}_{\hbar},
\end{equation} where the shearing places $\Mc$ in cohomological degree $-1$ and $\cO(\Mc)^{\shear}_{\hbar}$ is the Weyl algebra of $\Mc$ with shearing.

In particular, there is a tautological morphism \begin{equation}\label{eq:cohcorrloc}
\frc_{\Mc^*}^{\loc}: \Sat_{\hbar}(\Mc^*) * \delta_{X,\hbar}\to \delta_{X,\hbar}[1].
\end{equation}

\begin{remark}
\eqref{eq:cohcorrloc} is a special case of the \emph{local special cohomological correspondences} considered in \cite[\S5.1]{liu2025higherperiodintegralsderivatives}. While the description above determines $\frc_{\Mc^*}^{\loc}$ only up to a scalar, the fundamental class of the relative affine Grassmannian gives a canonical choice. Here we can choose $\frc_{\Mc^*}^{\loc}$ such that Proposition \ref{prop:momentloc} holds.
\end{remark}

\subsubsection{The Clifford algebra action}
Now we apply the construction of \cite{liu2025higherperiodintegralsderivatives} to obtain a Clifford algebra action on $\LL_G(\cP_X)$. Note that \cite{liu2025higherperiodintegralsderivatives} works with a sheaf theory that differs from ours in two respects: \begin{enumerate}
\item They work in the setting of $\ell$-adic sheaves, while we work with $D$-modules. This is not a problem, as all the constructions in \cite{liu2025higherperiodintegralsderivatives} can be carried out in the setting of ind-holonomic $D$-modules.
\item They use a version of $*$-sheaf theory (in which $*$-pullback is universally defined) instead of $!$-sheaf theory (in which $!$-pullback is universally defined). However, both sheaf theories define the same Plancherel algebra, which only depends on finite type geometry (see \cite[Remark\,8.3.2]{BZSV}). In fact, the only input from \cite{liu2025higherperiodintegralsderivatives} used in the article will be Proposition \ref{prop:clifaction} and Corollary \ref{cor:momentglob}, which involve only algebraic stacks on the automorphic side.
\end{enumerate}

Applying the local-to-global construction in \cite[\S4.4]{liu2025higherperiodintegralsderivatives} to \eqref{eq:cohcorrloc}, one obtains a morphism \begin{equation}
\label{eq:cohcorrglobC}
\frc_{\Mc^*}^{\glob}: \Sat_{C}(\Mc^*) * (\cP_X\boxtimes \omega_{C_{\dR}})\to \cP_X\boxtimes \omega_{C_{\dR}}[1]
.\end{equation}

Moreover, pushing forward along $C$ gives a morphism \begin{equation}
\label{eq:cohcorrglob}
\frc_{\Mc^*}: \Sat_{C}(\Mc^*) * \cP_X\to \cP_X[1]
.\end{equation} Translated to the spectral side, \eqref{eq:cohcorrglob} becomes \begin{equation}
\label{eq:cohcorrspec}
\frc_{\Mc^*}: \Gamma_{\dR}(C,\Mc_{\univ}^*)[-1]\otimes \LL_G(\cP_X)\to \LL_G(\cP_X)
.\end{equation}

Recall that $\Mc$ is a symplectic representation of $\Gc$. Therefore, the cohomology $\Gamma_{\dR}(C,\Mc_{\univ}^*)[-1]$ carries a symmetric bilinear form \[\cup:\Gamma_{\dR}(C,\Mc_{\univ}^*)^{\otimes 2} [-2]  \to \CC.\] Let \[\sw:\Gamma_{\dR}(C,\Mc_{\univ}^*)^{\otimes 2}\to \Gamma_{\dR}(C,\Mc_{\univ}^*)^{\otimes 2}\] be the swap map with the usual Koszul sign rule.
The following result is a direct translation of \cite[Corollary\,5.7]{liu2025higherperiodintegralsderivatives} to the spectral side:
\begin{prop}\label{prop:clifaction}
One has \[\frc_{\Mc^*}\circ \frc_{\Mc^*} - \frc_{\Mc^*}\circ \frc_{\Mc^*}\circ (\sw\otimes\id) = c \cdot \cup\otimes\id\in \Hom^0(\Gamma_{\dR}(C,\Mc_{\univ}^*)^{\otimes 2} [-2]\otimes \LL_G(\cP_X), \LL_G(\cP_X))\] where $c\in \CC^{\times}$ is a non-zero constant. In particular, the map $\frc_{\Mc^*}|_{\Loc_{\Gc}^{\circ}}$ induces an action of the Clifford algebra \[\mathbb{O}_{\Mc}^{\circ}:=\Cl(\Gamma_{\dR}(C,\Mc_{\univ}^*)[-1]|_{\Loc_{\Gc}^{\circ}})\] on $\LL_G(\cP_X)|_{\Loc_{\Gc}^{\circ}}$.
\end{prop}

\subsubsection{The moment map identity}
Consider the moment map $\mu_{\Mc}:\Mc\to \frgc^{*}$. It is a quadratic map, so it induces a morphism $\mu_{\Mc}^*:\frgc\to \Sym^2\Mc^* \sub (\Mc^*)^{\otimes 2}$. The following result can be directly deduced from \cite{Braverman_2022}:
\begin{prop}\label{prop:momentloc}
We have the following commutative diagram in $\Sph_{X,\hbar}^{\mathrm{nonren}}$:
\[\begin{tikzcd}
\Sat_{\hbar}(\frgc) * \delta_{X,\hbar} \ar[r, "\At^{\loc}_{\hbar}*\id"] \ar[d, "\Sat_{\hbar}(\mu_{\Mc}^*)*\id"] & \Sat_{\hbar}(\CC)[2] *\delta_{X,\hbar}\ar[d, "\sim"] \\
\Sat_{\hbar}((\Mc^*)^{\otimes 2}) * \delta_{X,\hbar} \ar[r, "\frc_{\Mc^*}^{\loc}\circ \frc_{\Mc^*}^{\loc}"] & \delta_{X,\hbar}[2]
\end{tikzcd}\] where $\frc_{\Mc^*}^{\loc}\circ \frc_{\Mc^*}^{\loc}: \Sat_{\hbar}((\Mc^*)^{\otimes 2}) * \delta_{X,\hbar}\to \delta_{X,\hbar}[2]$ is the map
\[\Sat_{\hbar}((\Mc^*)^{\otimes 2}) * \delta_{X,\hbar}\cong \Sat_{\hbar}(\Mc^*) * \Sat_{\hbar}(\Mc^*) * \delta_{X,\hbar}\xrightarrow{\id*\frc_{\Mc^*}^{\loc}}\Sat_{\hbar}(\Mc^*) * \delta_{X,\hbar}[1]\xrightarrow{\frc_{\Mc^*}^{\loc}}\delta_{X,\hbar}[2].\]
\end{prop}

\begin{proof}
After reduction modulo $\hbar$, this follows from \cite[Theorem\,2.2.1(a)]{Braverman_2022}. This implies the quantized version since \[\Hom^0_{\Sph_{X,\hbar}^{\mathrm{nonren}}}(\Sat_{\hbar}(\frgc)*\delta_{X,\hbar}, \delta_{X,\hbar})=0\] which can be seen from \eqref{eq:plXhbar}.

\end{proof}

Applying the local-to-global construction in \cite[\S4.4]{liu2025higherperiodintegralsderivatives}, we obtain the following global version of Proposition \ref{prop:momentloc}:
\begin{cor}\label{cor:momentglob}
We have the following commutative diagram in $\Dmod(\Bun_G\times C)$:
\[\begin{tikzcd}
\Sat_{C}(\frgc) * (\cP_X\boxtimes \omega_{C_{\dR}}) \ar[r, "\At^{\loc}_C*\id"] \ar[d, "\Sat_{C}(\mu_{\Mc}^*)*\id"] & \Sat_{C}(\CC)[2] * (\cP_X\boxtimes \omega_{C_{\dR}})\ar[d, "\sim"] \\
\Sat_{C}((\Mc^*)^{\otimes 2}) * (\cP_X\boxtimes \omega_{C_{\dR}}) \ar[r, "\frc_{\Mc^*}^{\glob}\circ \frc_{\Mc^*}^{\glob}"] & \cP_X\boxtimes \omega_{C_{\dR}}[2]
\end{tikzcd}.\] 
\end{cor}

This can be further translated to the spectral side, which gives the following result:
\begin{cor}\label{cor:momentglobspec}
We have the following commutative diagram in $\IndCoh_{\Nilp}(\Loc_{\Gc}\times C_{\dR})$:
\[\begin{tikzcd}
\frgc_{\univ} \otimes (\LL_G(\cP_X)\boxtimes \omega_{C_{\dR}}) \ar[r, "\At_C^{\loc}(\LL_G(\cP_X)\boxtimes \omega_{C_{\dR}})"] \ar[d, "(\mu_{\Mc}^*)_{\univ}\otimes\id"] & \LL_G(\cP_X)\boxtimes \omega_{C_{\dR}}[2] \ar[d, "\sim"] \\
(\Mc_{\univ}^*)^{\otimes 2} \otimes (\LL_G(\cP_X)\boxtimes \omega_{C_{\dR}}) \ar[r, "\frc_{\Mc^*}^{\glob}\circ \frc_{\Mc^*}^{\glob}"] & \LL_G(\cP_X)\boxtimes \omega_{C_{\dR}}[2]
\end{tikzcd}.\] 
\end{cor}

Taking cohomology along $C$ and applying Proposition \ref{prop:atloc=at}, we obtain the following result:
\begin{cor}\label{cor:momentglobint}
We have the following commutative diagram in $\IndCoh_{\Nilp}(\Loc_{\Gc})$:
\[\begin{tikzcd}
\Gamma_{\dR}(C,\frgc_{\univ})[-2] \otimes \LL_G(\cP_X) \ar[r, "\At(\LL_G(\cP_X))"] \ar[d, "{\Gamma_{\dR}(C,(\mu_{\Mc}^*)_{\univ})\otimes\id}"] & \LL_G(\cP_X) \ar[dd, "\sim"] \\
\Gamma_{\dR}(C, (\Mc_{\univ}^*)^{\otimes 2})[-2] \otimes \LL_G(\cP_X) \ar[d, "\Gys_{\D_C}\otimes\id "]  &  \\
\Gamma_{\dR}(C,\Mc^*_{\univ})^{\otimes 2}[-2]\otimes \LL_G(\cP_X)  \ar[r, "\frc_{\Mc^*}\circ \frc_{\Mc^*}"] & \LL_G(\cP_X)
\end{tikzcd}\] where $\Gys_{\D_C}:\Gamma_{\dR}(C, (\Mc_{\univ}^*)^{\otimes 2})\to \Gamma_{\dR}(C,\Mc^*_{\univ})^{\otimes 2}$ is the Gysin map induced by the $((\D_{C,\dR})_*,\D_{C,\dR}^!)$-adjunction, where $\D_C:C\to C\times C$ is the diagonal morphism.
\end{cor}

\subsubsection{Proof of Proposition \ref{prop:at=0}}\label{sec:proofat=0}
By Corollary \ref{cor:momentglobint}, it suffices to show that the morphism in $\QCoh(\Loc_{\Gc})$ \[\Gys_{\D_C}\circ \Gamma_{\dR}(C,(\mu_{\Mc}^*)_{\univ}): \Gamma_{\dR}(C,\frgc_{\univ})[-2] \to \Gamma_{\dR}(C,\Mc^*_{\univ})^{\otimes 2}[-2]\] is null-homotopic locally on $\Loc_{\Gc}^{\irr}$. This is obvious as the source is concentrated in cohomological degree $1$ while the target is concentrated in cohomological degree $0$ on the classical stack $\Loc_{\Gc}^{\irr}$.

\section{Classical limit}\label{sec:classicallimit}

In this section, we study the geometry of the relative Hitchin fibration, which serves as the microlocal picture of the period integral computation in \S\ref{sec:period}.

\subsection{The relative Hitchin fibration}
Let \[\Higgs_G:=\Map(C,[\frg^*_{\Omega_C}/G])\cong T^*\Bun_G\] be the Hitchin moduli stack, and let \[\Hitch_G:=\Map(C,\frg^*_{\Omega_C}/\!/ G)\] be the Hitchin base. Here the subscript $\Omega_C$ stands for twisting by the line bundle $\Omega_C$ via scalar dilation on $\frg^*$. Let \[h_G:\Higgs_G\to \Hitch_G\] be the Hitchin fibration.

\subsubsection{Geometry of the relative Hitchin fibration}

The relative Dolbeault moduli stack attached to the $G$-Hamiltonian space $M=T^*X$ is defined as
\[\Higgs_G^M = \Map(C, [M_{\Omega_C}/G])\cong \Bun_H^{(\frh')^*_{\Omega_C}}:=\Map(C,[(\frh')^*_{\Omega_C}/H])\] where $\frh' = \Lie(H')$ satisfies $\frh'|_H = \frh \oplus \Std\in \Rep(H)$ and $\Gm$ acts on $(\frh')^*$ by dilation.

We have a commutative diagram
\[\begin{tikzcd}
{[M/G]} \ar[r, "\mu_M"] \ar[d] & {[\frg^*/G]} \ar[dl] \\
M/\!/ G\cong \frg^*/\!/ G
\end{tikzcd}\] which induces a commutative diagram
\[\begin{tikzcd}
\Higgs_G^M \ar[d, "h_M"] \ar[r] & \Higgs_G \ar[dl, "h_G"] \\
\Hitch_G
\end{tikzcd}
.\] Here the isomorphism $M/\!/ G\cong \frg^*/\!/ G$ follows from the fact that $M$ is a \emph{strongly tempered} $G$-Hamiltonian space in the sense that $\Mc$ is a symplectic representation of $\Gc$. In this case, one has the following result:

\begin{prop}\label{prop:genfinite}
The map $h_M:\Higgs_G^M\to \Hitch_G$ is generically finite of degree $2^{(g-1)\dim(\Mc)}$.
\end{prop}
\begin{proof}
This follows from \cite[Corollary\,6.1]{hameister2024relative}. It can also be proved directly using the spectral curve description of $\Higgs_G^M$.
\end{proof}

\subsubsection{Coherence of pushforward}

Recall that we have fixed a square root $\Omega_C^{1/2}\in \Pic(C)$ of the canonical line bundle, which gives us a canonical square root $\cL_{G,\frg}^{1/2}\in \Pic(\Bun_G)$ of the determinant line bundle $\cL_{G,\frg}$ on $\Bun_G$ for any reductive group $G$. Since $\frh'|_H = \frh \oplus \Std\in \Rep(H)$, we get a canonical square root $\cL_{H,\Std}^{1/2}:=\cL_{H',\frh'}^{1/2}|_{\Bun_H}\otimes \cL_{H,\frh}^{-1/2}\in \Pic(\Bun_H)$ of the line bundle $\cL_{H,\Std}$. Via the map $a:\Higgs_G^M\cong \Bun_H^{(\frh')^*_{\Omega_C}}\to \Bun_H$, we regard $\cL_{H,\Std}^{1/2}$ as also a line bundle on $\Higgs_G^M$.

The main result of this section is the following theorem:

\begin{thm}\label{thm:classicallimcoh}
The object $(h_{M})_*(\cL_{H,\Std}^{-1/2})\in \QCoh(\Hitch_G)$ is a perfect complex.
\end{thm}

The idea of the proof of Theorem \ref{thm:classicallimcoh} is standard: using the Harder--Narasimhan stratification of $\Higgs_G^M$ (see \S\ref{sec:HNstrat}), we reduce to the semistable locus. Its good moduli space is proper over $\Hitch_G$, giving the required coherence. That is, consider the open substack $\Higgs_G^{M,\sst}:= \Bun_H^{(\frh')^*_{\Omega_C},\sst}\sub \Higgs_G^M$ in \S\ref{sec:HNstrat}. Let $h_M^{\sst}: \Higgs_G^{M,\sst}\to \Hitch_G$ be the restriction of $h_M$. We first prove the following:

\begin{lemma}\label{lem:ss=all}
The canonical map $(h_{M})_*(\cL_{H,\Std}^{-1/2})\to (h_{M}^{\sst})_*(\cL_{H,\Std}^{-1/2}|_{\Higgs_G^{M,\sst}})$ is an isomorphism.
\end{lemma}

Perfectness then follows from the following lemma, which is a special case of Proposition \ref{prop:HNstrat}, and the quasi-smoothness of $\Higgs_G^M$:
\begin{lemma}\label{lem:ssproper}
The following statements hold:
\begin{enumerate}
\item The underlying classical stack $\Higgs_G^{M,\sst,\cl}$ of $\Higgs_G^{M,\sst}$ admits a good moduli space $\phi:\Higgs_G^{M,\sst,\cl}\to \Higgs_G^{M,\sst,\cl,\coarse}$ in the sense of \cite[Definition\,4.1]{alper2013good}. In particular, $\phi_*$ preserves coherent sheaves.
\item The induced map $h_M^{\sst,\cl,\coarse}:\Higgs_G^{M,\sst,\cl,\coarse}\to \Hitch_G$ is proper.
\end{enumerate}
\end{lemma}

\subsubsection{A technical variant}

For the application in \S\ref{sec:period}, we also need a technical variant of Lemma \ref{lem:ss=all}. Let \begin{equation} \label{eq:bunbdd}j_{\bdd}:\Bun_H^{\bdd}\incl \Bun_H\end{equation} be the open substack where the Harder--Narasimhan filtration has slope $\lambda\in X_*(T_{H})_{\QQ}^+$ with $\lambda \leq (6g-6)\rho_H$ where $\rho_H$ is the half-sum of positive coroots of $H$ and $T_H$ is the maximal torus of $H$. Let $\Higgs_G^{M,\bdd}$ be its preimage under $\Bun_H^{(\frh')^*_{\Omega_C}}\to \Bun_H$. Let $h_M^{\bdd}:\Higgs_G^{M,\bdd}\to \Hitch_G$ be the restriction of $h_M$. We have the following lemma:

\begin{lemma}\label{lem:bdd=all}
The canonical map $(h_{M})_*(\cL_{H,\Std}^{-1/2})\to (h_{M}^{\bdd})_*(\cL_{H,\Std}^{-1/2}|_{\Higgs_G^{M,\bdd}})$ is an isomorphism.
\end{lemma}

We need the following variant of Theorem \ref{thm:classicallimcoh}, which is immediate from Theorem \ref{thm:classicallimcoh}, Lemma \ref{lem:bdd=all}, and Proposition \ref{prop:genfinite}:

\begin{cor}\label{cor:bddcoh}
The object $(h_{M}^{\bdd})_*(\cL_{H,\Std}^{-1/2}|_{\Higgs_G^{M,\bdd}}) \in \QCoh(\Hitch_G)$ is a perfect complex, generically concentrated in cohomological degree zero and of rank $2^{(g-1)\dim(\Mc)}$.
\end{cor}

The proofs of Lemma \ref{lem:bdd=all} and Lemma \ref{lem:ss=all} are in \S\ref{sec:bdd=all} and \S\ref{sec:ss=all}, respectively. They are direct applications of the method of \cite[Context\,6.5, Proposition\,6.6]{halpernleistner2023structuremoduligaugedmaps} using the Harder--Narasimhan stratification of $\Higgs_G^M$ (and $\Bun_H$) established in \S\ref{sec:HNstrat}.

\subsection{Harder--Narasimhan stratification}\label{sec:HNstrat}
In this section, we establish a Harder--Narasimhan stratification ($\Theta$-stratification) on moduli stacks of the form $\Bun_H^{(\frh')^*_{\Omega_C}}$. Although such a stratification appears to be well known to experts, we have not find a uniform treatment that applies directly to our setting.

Let $H$ be a connected semisimple group over $\CC$\footnote{Everything works for reductive groups over any field of characteristic zero with mild changes.} with a maximal torus $T$ and a Borel subgroup $B\supset T$, and let $V$ be a finite-dimensional representation of $H$. Let $\cL$ be a line bundle on $C$. We want to equip \[\Bun_H^{V_{\cL}}:=\Map(C,[V_{\cL}/H])\] with a stratification. Let \[\Hitch_H^{V_{\cL}} := \Map(C,V_{\cL}/\!/ H),\] and let \[h_H^{V_{\cL}}:\Bun_H^{V_{\cL}}\to \Hitch_H^{V_{\cL}}\] be the natural map.

Let $\lambda\in X_*(T)_{\QQ}^+$. For $U\in \Rep(H)$, we write $U^{\lambda \geq 0}$ for the sum of the weight spaces of non-negative weight with respect to $\lambda$, and we define $U^{\lambda > 0}$, $U^{\lambda = 0}$, $U^{\lambda \leq 0}$, $U^{\lambda <0}$ similarly. Let $P_{\lambda}$ be the parabolic subgroup of $H$ with $\frp_{\lambda}:=\Lie(P_{\lambda})=\frh^{\lambda \geq 0}$, and let $L_{\lambda}$ be its Levi subgroup such that $\frl_{\lambda}:=\Lie(L_{\lambda})=\frh^{\lambda =0}$. 
We have the following diagram:
\[\begin{tikzcd}
\Bun_{P_{\lambda}}^{V^{\lambda\geq 0}_{\cL}} \ar[r, "p_{\lambda}"] \ar[d, "q_{\lambda}"] & \Bun_H^{V_{\cL}} \\
\Bun_{L_{\lambda}}^{V^{\lambda=0}_{\cL}} \ar[u, bend left, "s_{\lambda}"] & 
\end{tikzcd}\]

There is a natural degree map $\deg:\Bun_{L_{\lambda}}\to X_*(Z(L_{\lambda}))_{\QQ}\sub X_*(T)_{\QQ}$. We use the same notation for the induced maps on $\Bun_{L_{\lambda}}^{V^{\lambda=0}_{\cL}}$ and $\Bun_{P_{\lambda}}^{V^{\lambda\geq 0}_{\cL}}$.

Let $(-,-)$ be the bilinear form on $X_*(T)_{\QQ}$ induced from the Killing form. For $\lambda\in X_*(T)_{\QQ}$, put $\lVert\lambda\rVert = (\lambda,\lambda)^{1/2}$.

Let $K/\CC$ be an algebraically closed field. For $y\in \Bun_{P_{\lambda}}^{V^{\lambda\geq 0}_{\cL}}(K)$, define
\[\mu(y) = \begin{cases} 0 & \text{if } \lambda = 0, \\[2pt] \dfrac{(\deg(y), \lambda)}{\lVert\lambda\rVert} & \text{otherwise}. \end{cases}\]
For $x\in \Bun_H^{V_{\cL}}(K)$, define
\[M(x) = \sup \{\mu(y)\mid \lambda\in X_*(T)_{\QQ}^+,\ y\in p_{\lambda}^{-1}(x)\}.\]
One has $M(x)\geq 0$ for every $x\in \Bun_H^{V_{\cL}}(K)$. If $M(x) = 0$, we say that $x$ is \emph{semistable}. More generally, for a reductive group $L$ (in our applications, a Levi subgroup of $H$) and an element $\theta\in X_*(Z(L))_{\QQ}$, replacing $\deg(y)$ by $\deg(y)-\theta$ in the definition of $\mu(y)$ defines $\theta$-semistability.

We say that $y\in \Bun_{P_{\lambda}}^{V^{\lambda\geq 0}_{\cL}}(K)$ is a \emph{Harder--Narasimhan filtration} of $x=p_{\lambda}(y)$ if $M(x) = \mu(y)=\lVert\lambda\rVert$.

The following result combines results from \cite{halpernleistner2023structuremoduligaugedmaps} and \cite{schmitt2008geometric}.

\begin{prop}\label{prop:HNstrat}
\leavevmode
\begin{enumerate}
\item For every $\lambda\in X_*(T)_{\QQ}^+$, the $\lambda$-semistable points form an open substack $\Bun_{L_\lambda}^{V_{\cL}^{\lambda = 0},\lambda\text{-}\sst}\subset \Bun_{L_\lambda}^{V_{\cL}^{\lambda = 0}}$ with quasi-compact connected components. In particular, the semistable points form a quasi-compact open substack $\Bun_H^{V_{\cL},\sst}\subset \Bun_H^{V_{\cL}}$.

\item Define $\Bun_H^{V_{\cL},(\lambda)} := q_{\lambda}^{-1}(\Bun_{L_{\lambda}}^{V_{\cL}^{\lambda = 0},\lambda\text{-}\sst})$. The map $p_{\lambda}:\Bun_{P_{\lambda}}^{V^{\lambda\geq 0}_{\cL}}\to \Bun_H^{V_{\cL}}$ restricts to a locally closed immersion $i_{\lambda}:\Bun_H^{V_{\cL},(\lambda)}\to \Bun_H^{V_{\cL}}$ on each connected component. These locally closed substacks form a stratification of $\Bun_H^{V_{\cL}}$. More precisely, the numerical invariant $\mu$ defines a $\Theta$-stratification in the sense of \cite[Theorem\,A]{halpernleistner2022structureinstabilitymodulitheory}.

\item The underlying classical stack $\Bun_H^{V_{\cL},\sst,\cl}$ of $\Bun_H^{V_{\cL},\sst}$ admits a good moduli space $\phi:\Bun_H^{V_{\cL},\sst,\cl}\to \Bun_H^{V_{\cL},\sst,\cl,\coarse}$ in the sense of \cite[Definition\,4.1]{alper2013good}. In particular, $\phi_*$ preserves coherent sheaves.

\item The induced map $h_H^{V_{\cL},\sst,\cl,\coarse}:\Bun_H^{V_{\cL},\sst,\cl,\coarse}\to\Hitch_H^{V_{\cL}}$ is proper.
\end{enumerate}
\end{prop}

\begin{remark}
Proposition \ref{prop:HNstrat} is an $\cL$-twisted version of the special case of \cite[Theorem\,1.2]{halpernleistner2023structuremoduligaugedmaps} where the $G$-scheme $X$ considered in \textit{loc. cit.} is a linear representation.
\end{remark}

\begin{proof}[Proof of Proposition \ref{prop:HNstrat}]
We only need to check that $\mu$ defines a $\Theta$-stratification. Given this, the quasi-compactness assertion in (1) and the good moduli space assertion in (3) follow from \cite[Theorem\,2.8.1.2]{schmitt2008geometric} and the fact that affine GIT quotients are good moduli spaces (see, for example, \cite[Remark\,4.8]{alper2013good}). (4) follows from \cite[Proposition\,2.8.1.4]{schmitt2008geometric}. The description of the strata using $\lambda$-semistable loci of the Levi subgroups follows from the same argument as in \cite[Proposition\,5.10]{halpernleistner2023structuremoduligaugedmaps}.

To prove that $\mu$ defines a $\Theta$-stratification, we need to check the conditions in \cite[Theorem\,A]{halpernleistner2022structureinstabilitymodulitheory}, and we will freely use the notations in \textit{loc. cit.} 

The condition $(\dagger\dagger)$ follows from that of $\Bun_H$ and the fact that $\Bun_H^{V_{\cL}}\to \Bun_H$ is representable and of finite type. The numerical invariant is defined from a line bundle by \cite[Lemma\,4.8]{halpernleistner2023structuremoduligaugedmaps} and hence is standard and satisfies condition $(R)$ by \cite[Lemma\,4.1.15, Lemma\,4.1.16]{halpernleistner2022structureinstabilitymodulitheory}.

For condition $(S)$, we show that for any discrete valuation ring $R/\CC$ with fraction field $K=\operatorname{Frac}(R)$ and residue field $\kappa$ and an $R$-point $x_R\in \Bun_H^{V_{\cL}}(R)$, given a Harder--Narasimhan filtration $y_K\in \Bun_{P_{\lambda}}^{V_{\cL}^{\lambda\geq 0}}(K)$ of the generic fiber $x_K$, there exists a filtration $y_{\kappa}\in \Bun_{P_{\lambda}}^{V_{\cL}^{\lambda\geq 0}}(\kappa)$ of the special fiber $x_{\kappa}$ such that \[\mu(y_{\kappa})\geq \mu(y_K).\] Moreover, when equality holds, the filtration extends to a filtration $y_{R}\in \Bun_{P_{\lambda}}^{V_{\cL}^{\lambda\geq 0}}(R)$. To prove this, we first consider the case where $V$ is trivial. Using the properness of the Drinfeld compactification, we can extend $y_K$ to a point $y_R'\in \overline{\Bun}_{P_{\lambda}}(R)$. Then we can take $y_{\kappa}$ to be the saturation of the special fiber $y_{\kappa}'\in \overline{\Bun}_{P_{\lambda}}(\kappa)$. This construction satisfies the desired properties. In the general case, since $V^{\lambda\geq 0}\sub V$ is a closed immersion, the condition on the $V$-section can be checked on the generic fiber, and hence, the same construction works. This proves condition $(S)$.

For condition $(B)$, we follow the strategy in the proof of \cite[Proposition\,4.6]{halpernleistner2023structuremoduligaugedmaps}. That is, we show that for any family $\xi:S\to \Bun_H^{V_{\cL}}$ for a quasi-compact scheme $S$, there exists a quasi-compact open substack $U\sub \Bun_H^{V_{\cL}}$ such that for any geometric point $s\in S(K)$ and a filtration $y\in \Bun_{P_{\lambda}}^{V_{\cL}^{\lambda\geq 0}}(K)$ of $\xi(s)$, there exists a filtration $y'\in\Bun_{P_{\lambda'}}^{V_{\cL}^{\lambda'\geq 0}}(K)$ of $\xi(s)$ with $\mu(y')\geq \mu(y)$ and $(p_{\lambda'}\circ s_{\lambda'}\circ q_{\lambda'})(y')\in U(K)$. When $V$ is trivial, this is \cite[Proposition\,4.6]{halpernleistner2023structuremoduligaugedmaps}. In general, we use a small variant of the argument. Indeed, we will choose $y'$ to be a coarsening of the filtration $y$ as follows. Take $y\in \Bun_{P_{\lambda}}^{V_{\cL}^{\lambda\geq 0}}(K)$ as above. If $\mu(y)\leq 0$, we can simply take $\lambda'=0$ and $y'$ to be the trivial filtration. So we can assume $\mu(y)>0$. Consider the hyperplanes in $X_*(T)_{\QQ}$ defined by the weights of $V$ and roots of $\frh$. They cut $X_*(T)_{\QQ}$ into cones. Let $\s_{\lambda}\sub X_*(T)_{\QQ}$ be the cone containing $\lambda$. Then for any coweight $\lambda'\in \overline{\s}_{\lambda}\backslash\{0\}$, the filtration $y$ induces a filtration $y'\in \Bun_{P_{\lambda'}}^{V_{\cL}^{\lambda'\geq 0}}(K)$ where $P_{\lambda}\sub P_{\lambda'}$ and $V^{\lambda\geq 0} \sub V^{\lambda'\geq 0}$. We choose $\lambda'\in \overline{\s}_{\lambda}\backslash\{0\}$ to maximize $\mu(y')=(\deg(y'),\lambda')/\lVert\lambda'\rVert=(\deg(y),\lambda')/\lVert\lambda'\rVert$. The maximizing property implies that we can rescale $\lambda'$ such that $(\lambda'-\deg(y),\lambda')=(\lambda'-\deg(y'),\lambda')=0$, where the last equality holds because $\deg(y')\in X_*(Z(L_{\lambda'}))_{\QQ}$ is the orthogonal projection of $\deg(y)\in X_*(Z(L_{\lambda}))_{\QQ}$ to $X_*(Z(L_{\lambda'}))_{\QQ}$. This implies $\lVert\lambda'\rVert = (\lambda',\lambda')/\lVert\lambda'\rVert=\mu(y')$. Since $S$ is quasi-compact, there exists a constant $C$ depending only on $S$, such that $\mu(y'')\leq C$ for any filtration $y''$ of the underlying $H$-bundle of $\xi(s)$ for any $s\in S(K)$. This implies that $\lVert\lambda'\rVert \leq C$. Moreover, since $\lambda'$ is the orthogonal projection of $\deg(y)$ on one of finitely many faces of the cones, we know that there exists an integer $N$ which only depends on $H, V$ such that $\lambda'\in \frac{1}{N}X_*(T)$. Therefore, we have $\lambda'\in \frac{1}{N}X_*(T)\cap B(0, C)$,\footnote{We write $B(x,r)$ for the ball of radius $r$ centered at $x$.} which is a finite set depending only on $S$; hence we can treat $\lambda'$ as fixed.
By \cite[Lemma\,4.10(2)]{halpernleistner2023structuremoduligaugedmaps}, we only need to show that there exists $d\in X_*(Z(L_{\lambda'}))_{\QQ}$ (depending only on $\lambda', C$) such that $\deg(y')-d\in X_*(Z(L_{\lambda'}))_{\QQ}^+$. Denote $\delta = \deg(y') - \lambda'$. Since $\lambda'$ is dominant in $X_*(Z(L_{\lambda'}))_{\QQ}$, we only need to show that there exists a constant $C'$ (depending only on $\lambda', C$) such that $\lVert\delta\rVert \leq C'$. Since $\lambda'$ is also strictly dominant in $X_*(Z(L_{\lambda'}))_{\QQ}$, there exists a real number $r>0$ such that $B(\lambda',r)\sub X_*(Z(L_{\lambda'}))_{\RR}^+$. By our choice of the constant $C$, we know that $(\deg(y'),-)\leq C \lVert-\rVert$ holds on $X_*(Z(L_{\lambda'}))_{\RR}^+$. Therefore, we have $(\delta,-) = (\deg(y'),-) - (\lambda',-)\leq (C+\lVert\lambda'\rVert)\lVert-\rVert\leq 2C\lVert-\rVert$ on $X_*(Z(L_{\lambda'}))_{\RR}^+$. In particular, we have $r\lVert\delta\rVert=(\delta,\lambda'+r\frac{\delta}{\lVert\delta\rVert})\leq 2C \lVert\lambda'+r\frac{\delta}{\lVert\delta\rVert}\rVert \leq 2C(C + r)$, which implies $\lVert\delta\rVert \leq 2C(C/r + 1)$. We can take $C' = 2C(C/r + 1)$ and conclude the proof.
\end{proof}

\subsection{Proof of Lemma \ref{lem:bdd=all}}\label{sec:bdd=all}
Consider the map $a:\Higgs_G^M \to \Bun_H$. Define
\[\cF = a_*(\cL_{H,\Std}^{-1/2}|_{\Higgs_G^M})\cong (a_*\cO_{\Higgs_G^M})\otimes \cL_{H,\Std}^{-1/2}.\]
We use the Harder--Narasimhan stratification of $\Bun_H$ as usual (same as taking $V$ to be the trivial representation in Proposition \ref{prop:HNstrat}). Let $i_{\lambda}:\Bun_H^{(\lambda)}\to \Bun_H$ be a Harder--Narasimhan stratum. Let $i_{\lambda}^{\wedge}:\Bun_H^{(\lambda),\wedge}\to \Bun_H$ be the formal completion of $\Bun_H$ along $\Bun_H^{(\lambda)}$. We have a diagram: \[
\begin{tikzcd}
\Bun_H^{(\lambda)} \ar[r, "i_{\lambda}"] \ar[d, "q_{\lambda}"] & \Bun_H \\
\Bun_{L_{\lambda}}^{\lambda-\sst}& 
\end{tikzcd}
\]

Let $X_*(T)_{\QQ}^{+,\bdd} \sub X_*(T)_{\QQ}^+$ be the subset of $\lambda\in X_*(T)_{\QQ}^+$ such that $\lambda\leq(6g-6)\rho_H$. Then we have $\Bun_H^{\bdd} = \bigsqcup_{\lambda\in X_*(T)_{\QQ}^{+,\bdd}} \Bun_H^{(\lambda)}$. For $\lambda\notin X_*(T)^{+,\bdd}_{\QQ}$, we want to show that $\Gamma(\Bun_H^{(\lambda),\wedge}, (i_{\lambda}^{\wedge})^!\cF) = 0$. For this, we only need to show that \[\Gamma(\Bun_H^{(\lambda)},i_{\lambda}^!\cF \otimes \Sym (T_{\Bun_H^{(\lambda)}/\Bun_H}[1]))=0.\]

Define $N=6g-6$. Since $\lambda\notin X_*(T)^{+,\bdd}_{\QQ}$, there exists a simple root $\alpha$ of $H$ such that $\langle \lambda, \alpha\rangle > N$. Choose $\varpi\in X_*(Z(L_{\lambda}))_{\QQ}$ via the following rule: Let $\varpi_1,\ldots,\varpi_r$ be the fundamental coweights of $\frh$ ordered according to the Dynkin diagram such that $\varpi_r$ (or $\varpi_{r-1}, \varpi_r$) corresponds to the spin node(s). Let $\alpha_1,\ldots,\alpha_r$ be the corresponding simple roots. Identifying $X_*(T)\cong \ZZ^r$ in the standard way, we have $\varpi_{r-1}= (1/2,\ldots,1/2,1/2,-1/2)$ and $\varpi_r = (1/2,\ldots,1/2,1/2,1/2)$ if $H$ is of type $D$. Let $\lambda = (a_1,\ldots,a_{r-1},a_r)$ with $a_1\geq \ldots \geq a_{r-1}\geq |a_r|$. Using the symmetry of the Dynkin diagram of type $D$, we can assume that $a_r\geq 0$.
\begin{itemize}
\item If $H$ is of type $B$, or $H$ is of type $D$ and $\alpha=\alpha_i\neq \alpha_{r-1},\alpha_r$, then let $\varpi = \varpi_i$.
\item If $H$ is of type $D$ and $a_{r-1}-a_r > N/3$, we take $\varpi = \varpi_{r-1} + \varpi_r$.
\item If $H$ is of type $D$ and $a_{r-1}-a_r \leq N/3$. In this case, one has $\alpha = \alpha_r$; therefore, $a_{r-1} + a_r > N$, which implies that $a_{r} > N/3$. We take $\varpi = 2\varpi_r$.
\end{itemize}

Consider the action of $\BB\Gm$ on $\Bun_{L_{\lambda}}^{\lambda-\sst}$ through a multiple of $\varpi$. We only need to show that $(q_{\lambda})_*(i_{\lambda}^!\cF \otimes \Sym (T_{\Bun_H^{(\lambda)}/\Bun_H}[1]))$ has (strictly) negative weights with respect to this $\BB\Gm$-action. 

Note that it admits a filtration with associated graded (up to a cohomological shift) \[\Sym(\Gamma(C,\frh_{\univ}\otimes \Omega_C)^*\oplus \Gamma(C,\Std_{\univ}\otimes \Omega_C)^*\oplus \Gamma(C,\frn_{\lambda,\univ}^{-})[1]\oplus \Gamma(C,\frn_{\lambda,\univ})^*[-1])\otimes \cL_{H,\Std}^{-1/2}\otimes \cL_{L_{\lambda},\frn_{\lambda}^-}^{-1}.\]
It remains to show that
\[\Sym(\Gamma(C,\frh_{\univ}^{\varpi<0}\otimes \Omega_C)^*\oplus \Gamma(C,\Std_{\univ}^{\varpi<0}\otimes \Omega_C)^*)\otimes \cL_{H,\Std}^{-1/2}\otimes \cL_{L_{\lambda},\frn_{\lambda}^-}^{-1}\] has strictly negative weights.

By our choice of $\varpi$, over $\Bun_{L_{\lambda}}^{\lambda-\sst}$, every fiber of $\frh^{\varpi<0}_{\univ}\otimes \Omega_C$ and $\Std^{\varpi<0}_{\univ}\otimes \Omega_C$ has strictly negative slopes. Therefore, the complex above equals \[
\Sym(H^1(C,\frh_{\univ}^{\varpi<0}\otimes \Omega_C)^*[1]\oplus H^1(C,\Std_{\univ}^{\varpi<0}\otimes \Omega_C)^*[1])\otimes \cL_{H,\Std}^{-1/2}\otimes \cL_{L_{\lambda},\frn_{\lambda}^-}^{-1}
.\]

Note that \[\max\weight_{\varpi}(\Sym(H^1(C,\frh_{\univ}^{\varpi<0}\otimes \Omega_C)^*[1])\otimes \cL_{L_{\lambda},\frn_{\lambda}^-}^{-1}) = \weight_{\varpi}(\cL_{L_{\lambda},\frn_{\lambda}}) - \weight_{\varpi}(\cL_{L_{\lambda},\frn_{\lambda}^-}) < 0\] and \[
\max\weight_{\varpi}(\Sym(H^1(C,\Std_{\univ}^{\varpi<0}\otimes \Omega_C)^*[1])\otimes \cL_{H,\Std}^{-1/2}) = (\weight_{\varpi}(\cL_{L_{\lambda},\Std^{\varpi>0}}) - \weight_{\varpi}(\cL_{L_{\lambda},\Std^{\varpi<0}}))/2 < 0
,\] where we use $\weight_{\varpi}$ to denote the weight with respect to the $\BB\Gm$-action and $\max\weight_{\varpi}$ for the maximal weight. This concludes the proof. \qed

\subsection{Proof of Lemma \ref{lem:ss=all}}\label{sec:ss=all}
Under the identification $\Higgs_G^M\cong \Bun_H^{(\frh')^*_{\Omega_C}}$ where $\frh' = \frh\oplus \Std$, we can use the Harder--Narasimhan stratification of $\Bun_H^{(\frh')^*_{\Omega_C}}$ as in Proposition \ref{prop:HNstrat}. Let $i_{\lambda}:\Bun_H^{(\frh')^*_{\Omega_C},(\lambda)}\to \Bun_H^{(\frh')^*_{\Omega_C}}$ be a Harder--Narasimhan stratum. We have a diagram:\[
\begin{tikzcd}
\Bun_H^{(\frh')^*_{\Omega_C},(\lambda)} \ar[r, "i_{\lambda}"] \ar[d, "q_{\lambda}"] & \Bun_H^{(\frh')^*_{\Omega_C}}  \\
\Bun_{L_{\lambda}}^{\frh'^{\lambda=0}_{\Omega_C},\lambda-\sst}  & 
\end{tikzcd}
\]

We proceed as in the proof of Lemma \ref{lem:bdd=all}. We only need to show that for any nonzero $\lambda\in X_*(T)_{\QQ}^+$ \[\Gamma(\Bun_H^{(\frh')^*_{\Omega_C},(\lambda)}, i_{\lambda}^!\cL_{H,\Std}^{-1/2}\otimes \Sym(T_{\Bun_H^{(\frh')^*_{\Omega_C},(\lambda)}/\Bun_H^{(\frh')^*_{\Omega_C}}}[1]))=0.\]

Consider the action of $\BB\Gm$ on $\Bun_{L_{\lambda}}^{\frh'^{\lambda=0}_{\Omega_C},\lambda-\sst}$ through a multiple of $\lambda$. We only need to show that \[(q_{\lambda})_*(i_{\lambda}^!\cL_{H,\Std}^{-1/2}\otimes \Sym(T_{\Bun_H^{(\frh')^*_{\Omega_C},(\lambda)}/\Bun_H^{(\frh')^*_{\Omega_C}}}[1]))\] has strictly negative weights with respect to this $\BB\Gm$-action. This complex has a filtration with associated graded (up to a cohomological shift) \[\Sym(\Gamma(C,[\frn_{\lambda,\univ}^-\to (\frh'_{\univ})^{\lambda<0}\otimes\Omega_C])\oplus \Gamma(C,[\frn_{\lambda,\univ}\to (\frh'_{\univ})^{\lambda>0}\otimes\Omega_C])^*
)\otimes \cL_{H,\Std}^{-1/2}\otimes \cL_{L_{\lambda},\frn_{\lambda}^-}^{-1}\otimes \cL_{L_{\lambda},\frh'^{\lambda<0}_{\Omega_C}}\] where both $[\frn_{\lambda,\univ}^-\to (\frh'_{\univ})^{\lambda<0}\otimes\Omega_C]$ and $[\frn_{\lambda,\univ}\to (\frh'_{\univ})^{\lambda>0}\otimes\Omega_C]$ are placed in degrees $-1, 0$. The desired vanishing follows because the maximal weight is \[\weight_{\lambda}(\cL_{H,\Std}^{-1/2}\otimes \cL_{L_{\lambda},\frn_{\lambda}^-}^{-1}\otimes \cL_{L_{\lambda},\frh'^{\lambda<0}_{\Omega_C}})=(\weight_{\lambda}(\cL_{L_{\lambda},\Std^{\lambda>0}})-\weight_{\lambda}(\cL_{L_{\lambda},\Std^{\lambda<0}}))/2 + (\weight_{\lambda}(\cL_{L_{\lambda},\frn_{\lambda}}) - \weight_{\lambda}(\cL_{L_{\lambda},\frn_{\lambda}^-})) < 0.\] This concludes the proof. \qed

\section{Period of differential operators}\label{sec:period}

In this section, we compute the $*$-period integral of ${\rm D}_{\Bun_G,\kappa}\in\Dmod(\Bun_G)$. More precisely, let $\Op_{\Gc}$ be the moduli space of $\Gc$-opers on $C$ and $p_{\Op}:\Op_{\Gc}\to \Loc_{\Gc}$ be the natural map.\footnote{We use the standard convention for opers as in \cite[\S10.2]{gaitsgory2015outline} which is different from the one used in \cite{beilinson1991quantization} when $\Gc$ is not of adjoint type.} We will prove the following:

\begin{prop}\label{prop:onefiber}
There exists $\s\in \Op_{\Gc}(\CC)$ such that the fiber $i_{\s}^*(\LL_G(\cP_X))[\delta_H-\delta_{N_{\Omega_C^{1/2}}}]\in \Vect$ lies in $\Vect^{\heartsuit}$ and has dimension $2^{(g-1)\dim(\Mc)}$, where $i_{\s}:\{\s\}\to \Loc_{\Gc}$ is the induced map.
\end{prop}

We prove Proposition \ref{prop:onefiber} by passing to the classical limit, where the corresponding statement is Corollary \ref{cor:bddcoh}.

\subsection{Filtered $D$-modules}
We first recall the theory of filtered $D$-modules on smooth algebraic stacks following \cite[\S6.2]{faergeman2025non}.

\subsubsection{Hodge stack}\label{sec:hodge}

For a smooth algebraic stack $\mathscr{Y}$, there is a \emph{Hodge stack} ${\scr Y}_{\hbar,\dR}\to [\AA^1_{\hbar}/\Gm]$. Its fiber over $\hbar=1$ is $\mathscr{Y}_{\dR}$, and its fiber over $\hbar=0$ is $\mathbb{B}_{\scr Y}T_{\scr Y}^{\wedge}$. The category of filtered $D$-modules is defined by \[\Dmod_{\hbar}(\mathscr{Y}) := \IndCoh(\mathscr{Y}_{\hbar,\dR}).\] 

We have natural maps $j_1:{\scr Y}_{\dR}\to {\scr Y}_{\hbar,\dR}$ and $i_0:\mathbb{B}_{\scr Y}T_{\scr Y}^{\wedge}\to {\scr Y}_{\hbar,\dR}$. We have the forgetful functor \[\oblv:=j_1^!:\Dmod_{\hbar}(\mathscr{Y})\cong \IndCoh({\scr Y}_{\hbar,\dR}) \to \Dmod(\mathscr{Y})\cong \IndCoh({\scr Y}_{\dR}).\] Koszul duality gives $\IndCoh(\mathbb{B}_{\scr Y}T_{\scr Y}^{\wedge})\cong \QCoh(T^*{\scr Y})$. We have the associated graded functor \[\gr:=i_0^!:\Dmod_{\hbar}(\mathscr{Y})\cong \IndCoh({\scr Y}_{\hbar,\dR}) \to \QCoh(T^*{\scr Y})\cong \IndCoh(\mathbb{B}_{\scr Y}T_{\scr Y}^{\wedge}).\]

\subsubsection{Functoriality}\label{sec:filfunc}

Given a morphism between smooth algebraic stacks $f:{\scr Y} \to {\scr Z}$, we have the induced correspondence on cotangent stacks \[\begin{tikzcd}
T^*{\scr Z}\times_{\scr Z}{\scr Y} \ar[r, "Pf"] \ar[d, "Df"] & T^*{\scr Z} \\
T^*{\scr Y} & 
\end{tikzcd}.\]

We use \[f_{\dR,*,\ren}:\Dmod({\scr Y}) \to \Dmod({\scr Z})\] to denote the renormalized direct image functor defined by continuously extending the functor \[f_{\dR,*}:\Dmod({\scr Y})\to \Dmod({\scr Z})\] on compact objects. The functors \[f_{\dR}^!:\Dmod({\scr Z}) \to \Dmod({\scr Y})\] and \[f_{\dR,*,\ren}:\Dmod({\scr Y}) \to \Dmod({\scr Z})\] admit natural filtered enhancements, for which we use the same notation. By \cite[(6.2.1)]{faergeman2025non}, one has commutative diagrams of functors:
\[\begin{tikzcd}
\Dmod_{\hbar}({\scr Z}) \ar[r, "f_{\dR}^!"] \ar[d, "\gr"] & \Dmod_{\hbar}({\scr Y}) \ar[d, "\gr"] \\
\QCoh(T^*{\scr Z}) \ar[r, "(Df)_*(Pf)^!"] & \QCoh(T^*{\scr Y})
\end{tikzcd} \]
\[\begin{tikzcd}
\Dmod_{\hbar}({\scr Y}) \ar[r, "f_{\dR, *,\ren}"] \ar[d, "\gr"] & \Dmod_{\hbar}({\scr Z}) \ar[d, "\gr"] \\
\QCoh(T^*{\scr Y}) \ar[r, "(Pf)_*(Df)^*"] & \QCoh(T^*{\scr Z})
\end{tikzcd} .\]

\subsubsection{Filtration on differential operators}\label{sec:fildif}

Let ${\scr Y}$ be a smooth algebraic stack. There is a natural map $q_{{\scr Y},\hbar}:{\scr Y}\times [\AA^1_{\hbar}/\Gm] \to {\scr Y}_{\hbar,\dR}$. This induces an adjoint pair \[\ind_{\hbar}^r:=(q_{{\scr Y},\hbar})_*: \IndCoh({\scr Y}\times [\AA^1_{\hbar}/\Gm]) \rightleftarrows \Dmod_{\hbar}({\scr Y}) : \oblv_{\hbar}^r:=(q_{{\scr Y, \hbar}})^!.\] 
In particular, we have the filtered $D$-module ${\rm D}_{\hbar,\scr Y}:=\ind_{\hbar}^r(\omega_{{\scr Y}\times [\AA^1_{\hbar}/\Gm]})\in \Dmod_{\hbar}({\scr Y})$ satisfying \[\oblv({\rm D}_{\hbar,\scr Y})\cong {\rm D}_{\scr Y}\in \Dmod({\scr Y})\] and \[\gr({\rm D}_{\hbar,\scr Y})\cong \pi_{{\scr Y}}^*\omega_{\scr Y}\] where $\pi_{\scr Y}:T^*{\scr Y}\to {\scr Y}$ is the projection. More generally, for any line bundle $\cL\in \Pic({\scr Y})$, one can define the category of twisted filtered $D$-modules $\Dmod_{\hbar,\cL}({\scr Y})$ and the object ${\rm D}_{\hbar,{\scr Y},\cL}\in \Dmod_{\hbar,\cL}({\scr Y})$ via the standard twisting construction. Using the canonical equivalence $\Dmod_{\hbar,\cL}({\scr Y}) \cong \Dmod_{\hbar}({\scr Y})$, we have \[\oblv({\rm D}_{\hbar,{\scr Y},\cL}) \cong {\rm D}_{{\scr Y},\cL}\in \Dmod_{\cL}({\scr Y})\] and \[\gr({\rm D}_{\hbar,{\scr Y},\cL}) \cong \pi_{\scr Y}^*(\cL^{-1}\otimes \omega_{\scr Y}).\]

\subsection{Quantization of the Hitchin system}\label{sec:quanthit}
In this section, we review the result of \cite{beilinson1991quantization} needed for the computation.

The following proposition follows from \cite[Theorem\,18.5.2]{arinkin2024proof}.\footnote{We ignore some line twists in \cite[Theorem\,18.5.2]{arinkin2024proof} which are not relevant for our purpose.}
\begin{prop}\label{prop:dbungop}
There is an isomorphism $\LL_{G}({\rm D}_{\Bun_G, \kappa})\cong (p_{\Op})_*\cO_{\Op_{\Gc}}[\delta_{N_{\Omega_C^{1/2}}}]\in \IndCoh_{\Nilp}(\Loc_{\Gc})$.
\end{prop}

In particular, there is a natural map \[\cO(\Op_{\Gc}) \to \End({\rm D}_{\Bun_G, \kappa}).\] Recall that we have the Hitchin map $h_G:\Higgs_G\to \Hitch_G$. The main result of \cite{beilinson1991quantization} says that $\cO(\Op_{\Gc}) \to \End({\rm D}_{\Bun_G, \kappa})$ is a quantization of the map \[\cO(\Hitch_G) \to \End((\cL_{G,\frg}^{-1/2}\otimes\omega_{\Bun_G})|_{\Higgs_G})\] induced by $h_G$. 

More precisely, let $\Op_{\Gc,\hbar}$ be the moduli of $\hbar$-opers for $\Gc$ (replacing connections by $\hbar$-connections), which is a family over $[\AA^1_{\hbar}/\Gm]$, whose fiber at $\hbar=0$ is $\Hitch_G$ and the fiber at $\hbar=1$ is $\Op_{\Gc}$. Let ${\rm D}_{\hbar,\Bun_G,\kappa}\in \Dmod_{\hbar}(\Bun_G)$ be the natural filtered $D$-module on $\Bun_G$, whose fiber at $\hbar=0$ is $(\cL_{G,\frg}^{-1/2}\otimes\omega_{\Bun_G})|_{\Higgs_G}\in \QCoh(\Higgs_G)$ and the fiber at $\hbar=1$ is ${\rm D}_{\Bun_G,\kappa}\in \Dmod(\Bun_G)$. The following is a reformulation of the main result of \cite{beilinson1991quantization}.

\begin{prop}\label{prop:bdquant}
There is a natural map of algebra objects $\cO(\Op_{\Gc,\hbar}) \to \underline{\End}_{\QCoh([\AA^1_{\hbar}/\Gm])}({\rm D}_{\hbar,\Bun_G, \kappa})$, whose fibers at $\hbar=0$ and $\hbar=1$ are the maps considered above.
\end{prop}

\subsection{Computation of the period of differential operators}
In this section, we compute the $*$-period integral of ${\rm D}_{\Bun_G,\kappa}\in \Dmod(\Bun_G)$. More precisely, we will prove the following result.
\begin{prop}\label{prop:periodcomp}
The object \[\Hom(\cP_X,{\rm D}_{\Bun_G,\kappa})\in \Mod_{\cO(\Op_{\Gc})}\] admits an increasing filtration that is bounded below such that the associated graded is isomorphic to \[\Gamma(\Higgs_G^{M,\bdd},\cL_{H,\Std}^{-1/2}|_{\Higgs_G^{M,\bdd}})[\delta_H]\in \Mod_{\cO(\Hitch_G)}.\]
\end{prop}

\begin{proof}

We start with the following lemma. Recall $j_{\cusp}:\Bun_G^{\cusp}\to \Bun_G$ defined in \eqref{eq:buncusp}.

\begin{lemma}\label{lem:dbuncusp}

The object ${\rm D}_{\Bun_G,\kappa}\in \Dmod(\Bun_G)$ is cuspidal. In particular, one has \[{\rm D}_{\Bun_G,\kappa} \simeq (j_{\cusp})_{*}{\rm D}_{\Bun_G^{\cusp},\kappa}.\] 
\end{lemma}
\begin{proof}
Since the map $p_{\Op}:\Op_{\Gc}\to \Loc_{\Gc}$ factors through $\Loc_{\Gc}^{\irr}\sub \Loc_{\Gc}$, the statement follows from \cite[Corollary\,11.2.2.3]{faergeman2025non}.
\end{proof}

Recall $j_{\bdd}:\Bun_H^{\bdd}\to \Bun_H$ defined in \eqref{eq:bunbdd}. The restriction of $\pi_X:\Bun_H\to \Bun_G$ induces a map $\pi_X^{\bdd}:\Bun_H^{\bdd}\to \Bun_G$ which satisfies $\pi_X^{-1}(\Bun_G^{\cusp}) \sub \Bun_H^{\bdd}$. By Lemma \ref{lem:dbuncusp}, we have \begin{equation}\label{eq:computation1}\begin{split}
&\Hom(\cP_X,{\rm D}_{\Bun_G,\kappa}) \\ 
 \cong & \Gamma_{\dR}(\Bun_G^X,\pi_X^!{\rm D}_{\Bun_G,\kappa}) \\
 \cong & \Gamma_{\dR}(\Bun_H^{\bdd}, (\pi_X^{\bdd})^!{\rm D}_{\Bun_G,\kappa})
\end{split}.\end{equation}
We note the following lemma.

\begin{lemma}\label{lem:dr=ren}
We have $\Gamma_{\dR,\ren}(\Bun_H^{\bdd}, (\pi_X^{\bdd})^!{\rm D}_{\Bun_G,\kappa}) \cong \Gamma_{\dR}(\Bun_H^{\bdd}, (\pi_X^{\bdd})^!{\rm D}_{\Bun_G,\kappa})$.
\end{lemma}
\begin{proof}
We note the following general fact: Let $\mathscr Y$ be a smooth quasi-compact algebraic stack with affine diagonal. By \cite[Theorem\,1.4.2,\,Theorem\,3.3.5]{drinfeld2013some}, the category $\QCoh(\mathscr Y)$ is compactly generated and $\Gamma(\mathscr Y,-):\QCoh(\mathscr Y)\to\Vect$ is continuous. Therefore, we have \[\Gamma_{\dR,\ren}({\scr Y},-)\circ \ind^l \cong \Gamma({\scr Y}, -\otimes \omega_{\scr Y})\cong  \Gamma_{\dR}({\scr Y},-)\circ \ind^l.\] Since $\Bun_H^{\bdd}$ is QCA, we only need to show that \[(\pi_X^{\bdd})^!{\rm D}_{\Bun_G,\kappa}\in \ind^l(\QCoh(\Bun_H^{\bdd}))\sub \Dmod(\Bun_H^{\bdd}).\] 
Let $f_{\bdd}:\Bun_H^{\bdd}\to \Bun_{H'}$ be the natural map. We have
\[\begin{split}
& (\pi_X^{\bdd})^!{\rm D}_{\Bun_G,\kappa} \\
\cong & (\pi_X^{\bdd})^!\ind^l(\cL_{G,\frg}^{-1/2}) \\
\cong & \ind^l(\cL_{H,\frh}^{-1/2}|_{\Bun_H^{\bdd}})\otimes^! f_{\bdd}^!\ind^l(\cL_{H',\frh'}^{-1/2}) \\
\cong & \ind^l(\cL_{H,\frh}^{-1/2}|_{\Bun_H^{\bdd}}\otimes \oblv^lf_{\bdd}^!\ind^l(\cL_{H',\frh'}^{-1/2})) \\
\end{split}.
\]
This concludes the proof.
\end{proof}

By the discussion in \S\ref{sec:filfunc}, we obtain a natural increasing filtration on $\Gamma_{\dR,\ren}(\Bun_H^{\bdd}, (\pi_X^{\bdd})^!{\rm D}_{\Bun_G,\kappa})$. By \eqref{eq:computation1} and Lemma \ref{lem:dr=ren}, we obtain the natural increasing filtration on $\Hom(\cP_X,{\rm D}_{\Bun_G,\kappa})$ compatible with the filtration on $\cO(\Op_{\Gc})$. The filtration is bounded below by construction. Consider the diagram of cotangent stacks

\[\begin{tikzcd}
\Higgs_G^{M,\bdd} \ar[r] \ar[d] & T^*\Bun_G\times_{\Bun_G} \Bun_H^{\bdd} \ar[r, "P\pi_X^{\bdd}"] \ar[d, "D\pi_X^{\bdd}"] & T^*\Bun_G  \\
\Bun_H^{\bdd} \ar[r ,"Dp_X"] \ar[d, "Pp_X"] & T^*\Bun_H^{\bdd} & \\
* && 
\end{tikzcd}
\] where the square is Cartesian and $p_X:\Bun_H^{\bdd}\to *$ is the map to a point. By the discussion in \S\ref{sec:filfunc} and \S\ref{sec:fildif}, we have an isomorphism

\[
\begin{split}
& \gr(\Hom(\cP_X,{\rm D}_{\Bun_G,\kappa})) \\
\cong & \Gamma(\Bun_H^{\bdd}, (Dp_X)^*(D\pi_X^{\bdd})_*(P\pi_X^{\bdd})^!(\cL_{G,\frg}^{-1/2}\otimes \omega_{\Bun_G}|_{T^*\Bun_G})) \\
\cong & \Gamma(\Higgs_G^{M,\bdd}, \cL_{G,\frg}^{-1/2}\otimes \omega_{\Bun_H^{\bdd}}|_{\Higgs_G^{M,\bdd}}) \\
\cong & \Gamma(\Higgs_G^{M,\bdd}, \cL_{H,\Std}^{-1/2}|_{\Higgs_G^{M,\bdd}})[\delta_H]\\
\end{split}
\] which is compatible with the $\cO(\Hitch_G)$-module structure by Proposition \ref{prop:bdquant}. This concludes the proof.

\end{proof}

\begin{remark}\label{rmk:comparefr}
    Our proof of Proposition \ref{prop:periodcomp} can be compared with the proof of \cite[Theorem\,6.1.2.1]{faergeman2025non} for the Whittaker period. Both proofs pass to the classical limit, but the technical subtleties are different. In our setting, the filtration is automatically bounded below, unlike in \cite{faergeman2025non}, where the microlocal geometry is non-polarized. On the other hand, we need to replace $\Bun_H$ with a quasi-compact stack $\Bun_H^{\bdd}$ and establish Lemma \ref{lem:dr=ren} to ensure the exhaustiveness of the filtration, whereas exhaustiveness is automatic in \cite{faergeman2025non}. Moreover, the classical limit is considerably more complicated in our setting. Indeed, the relative Hitchin fibration is much more involved than in the Whittaker case, where the corresponding $\Higgs_G^M$ is the Kostant section. This leads to the additional calculations in \S\ref{sec:classicallimit}.
\end{remark}

\subsection{Proof of Proposition \ref{prop:onefiber}}
The following lemma is standard:
\begin{lemma}\label{lem:perfectmodhbar}
Let $A\in \operatorname{Alg}(\QCoh([\AA^1_{\hbar}/\Gm]))$ be a filtered algebra, and let $M\in \Mod_A$ be a filtered $A$-module. Assume that:
\begin{enumerate}
\item The filtrations on $A$ and $M$ are increasing and bounded below.
\item $A/\hbar A \in \Vect^{\heartsuit}$.
\item The associated graded $M/\hbar M \in \Mod_{A/\hbar A}$ is perfect.
\end{enumerate}
Then $M$ is perfect as an $A$-module.
\end{lemma}
\begin{proof}
Note that $A/\hbar A \in \Vect^{\heartsuit}$ implies that $H^i(A)/\hbar H^i(A)=0$ for $i\neq 0$, which implies that $H^i(A)=0$ for $i\neq 0$. Therefore, we have $A\in \Vect^{\heartsuit}$. Assume $M/\hbar M$ has Tor amplitude $[a,b]$. By the same argument, we know that $M\in\Mod_A^{\leq b}$. Take a finitely generated projective module $\overline{P}\in \Mod_{A/\hbar A}^{\heartsuit}$ and a map $\overline{P}[-b]\to M/\hbar M$ inducing $\overline{P}\surj H^b(M/\hbar M)$  (when $a=b$, we simply take $\overline{P} = (M/\hbar M)[b]$). We can find a finitely generated projective module $P\in \Mod_A^{\heartsuit}$ and a lift $P[-b]\to M$. Then $\operatorname{Cone}(P[-b]\to M)/\hbar$ has Tor amplitude $[a,b-1]$. Iterating this process shows that $M$ is perfect.
\end{proof}

By Lemma \ref{lem:perfectmodhbar} and Corollary \ref{cor:bddcoh}, we know that $\Hom(\cP_X, {\rm D}_{\Bun_G,\kappa})\in \Mod_{\cO(\Op_{\Gc})}$ is perfect. Moreover, there exists $\s\in \Op_{\Gc}(\CC)$ such that $i_{\s}^*(\Hom(\cP_X, {\rm D}_{\Bun_G,\kappa}))[-\delta_H]\in \Vect^{\heartsuit}$ and has dimension $2^{(g-1)\dim(\Mc)}$. We have \[\Hom(\cP_X, {\rm D}_{\Bun_G,\kappa})\cong \Hom(\LL_G(\cP_X),(p_{\Op})_*\cO_{\Op_{\Gc}}[\delta_{N_{\Omega_C^{1/2}}}])\cong (p_{\Op}^*\LL_{G}(\cP_X))^{\vee}[\delta_{N_{\Omega_C^{1/2}}}]\in \QCoh(\Op_{\Gc})\] where $(p_{\Op}^*\LL_{G}(\cP_X))^{\vee}$ is the linear dual of $p_{\Op}^*\LL_{G}(\cP_X)\in \QCoh(\Op_{\Gc})$ which is well-defined by Lemma \ref{lem:coherentcoh}. This implies that \[i_{\s}^*(\LL_G(\cP_X))[\delta_H-\delta_{N_{\Omega_C^{1/2}}}]\cong i_{\s}^*(\Hom(\cP_X, {\rm D}_{\Bun_G,\kappa}))^{\vee}[\delta_H]\in \Vect^{\heartsuit}\] and has dimension $2^{(g-1)\dim(\Mc)}$. This concludes the proof of Proposition \ref{prop:onefiber}. \qed

\bibliographystyle{amsalpha}
\bibliography{Bibliography}
\end{document}